\documentclass[12pt]{amsart}  
\usepackage{amsmath,amsthm,amsfonts,amssymb, eucal}
\numberwithin{equation}{section}
\newtheorem{thm}[equation]{Theorem}
\newtheorem{lem}[equation]{Lemma}
\newtheorem{prop}[equation]{Proposition}
\newtheorem{ex}[equation]{Example}
\newtheorem{exs}[equation]{Examples}
\theoremstyle{definition} 
\newtheorem{rem}[equation]{Remark}

\newcommand \coop {\text{coop}}
\newcommand{\e}{\epsilon}
\newcommand{\E}{\mathcal E}
\newcommand{\F}{\mathcal F}

\newcommand{\id}{\hbox{\rm id}}
\newcommand{\kc}{\mathbb {K}}
\newcommand{\la}{\langle}
\newcommand{\m} {\medskip}
\newcommand \op {\text{op}}
\newcommand{\ot}{\otimes}
\newcommand{\ra}{\rangle}
\newcommand{\w} {\omega}
\newcommand{\K}{\mathbb {K}}
\newcommand{\CN}{\mathcal{N}}

\newcommand{\KK}{\mathfrak {K}}
\newcommand{\UK}{U_{\KK}}
\newcommand{\Z}{\mathbb {Z}}

\begin{document}

\title[DEFORMATIONS OF ALGEBRAS]
{QUANTUM GROUP ACTIONS,  \\ TWISTING ELEMENTS, and \\
DEFORMATIONS OF ALGEBRAS}

\author{Georgia Benkart} \address{Department of Mathematics \\ University
of Wisconsin \\  Madison, Wisconsin 53706, USA}
\email{benkart@math.wisc.edu}  

\author{Sarah Witherspoon}
\address{Department of Mathematics\\
 Texas A\&M University \\ College Station, Texas 77843, USA}
\email{sjw@math.tamu.edu}
\thanks{The authors gratefully acknowledge support from the following
grants:\quad NSF grants \#{}DMS--0245082 and \#{}DMS--0443476, and
NSA grants MSPF-02G-082 and
MDA904-03-1-0068. The second author thanks the Mathematical Sciences Research
Institute for its support during the writing of this paper.}

\date{December 6, 2005}


\begin{abstract}
We construct twisting elements for module algebras
of restricted two-parameter quantum groups from factors of their
$R$-matrices. We generalize the theory of Giaquinto and Zhang to
universal deformation formulas for categories of module algebras
and give examples arising from $R$-matrices of two-parameter quantum
groups.
\end{abstract}

\maketitle

\section{Introduction}    

The algebraic deformation theory
of Giaquinto and Zhang  \cite{giaquinto-zhang98}   laid out a
framework for certain types of deformations of associative algebras,
namely those arising from actions of bialgebras on algebras.
This work generalized known
examples of deformations coming from the action of a Lie algebra (equivalently of its
universal enveloping algebra).
Giaquinto and Zhang illustrated their ideas by constructing a 
new deformation
formula from the action of a noncommutative bialgebra, namely, the universal
enveloping algebra of a certain nonabelian Lie algebra.  This
Lie algebra is spanned by the $n \times n$ matrix units
$E_{1,p},\ E_{p,n}$ for $p =1, \dots, n$,  along with the $n \times n$
diagonal matrices of trace 0, and so it is an
abelian extension of a Heisenberg Lie algebra when $n \geq 3$,  and it is
a two-dimensional nonabelian Lie algebra when $n =2$.
In recent work, Grunspan \cite{grunspan} applied the deformation formula
for the $n=2$ case  to solve the open problem of giving  an explicit
deformation of the Witt algebra.

In their study of deformations and orbifolds,
C\u ald\u araru, Giaquinto, and the second author
\cite{caldararu-giaquinto-witherspoon}
constructed a deformation from a non{\em co}commutative
(and noncommutative) bialgebra -- the Drinfeld double of a Taft algebra.
This  was related to unpublished work of Giaquinto and Zhang and
was the first known example of an explicit formula for a formal
deformation of an algebra arising from a noncocommutative bialgebra.

The present paper began as an attempt to put the example of
\cite{caldararu-giaquinto-witherspoon} into a general context.
The algebra in \cite{caldararu-giaquinto-witherspoon}  to undergo deformation
was a crossed
product of a polynomial algebra and the group algebra of a finite group. We
expected that  
generalizations of the Drinfeld double of a Taft algebra, such as the
two-parameter restricted quantum groups ${\mathfrak{u}}_{r,s}
({\mathfrak{sl}}_{n})$
of our paper \cite{benkart-witherspoon3},
would act on similar crossed products, potentially leading to
deformations. Indeed, ${\mathfrak{u}}_{r,s}({\mathfrak{sl}}_{n})$
{\em does}
act on such a crossed product, but we were unable to generalize
the deformation formula of \cite{caldararu-giaquinto-witherspoon}
directly to deform the multiplication in the crossed product  via the
${\mathfrak{u}}_{r,s}({\mathfrak{sl}}_{n})$-action.
Instead, a generalization in a different direction proved more
successful in deforming these particular types of algebras (see
\cite{witherspoon}).
However, as we discuss in Sections 3 through 5 in this paper,  certain
related noncommutative algebras and their crossed products with group algebras
also carry a ${\mathfrak{u}}_{r,s}({\mathfrak{sl}}_n)$-action.  
Since ${\mathfrak{u}}_{r,s}({\mathfrak{sl}}_n)$ is a quotient of
the infinite-dimensional two-parameter quantum group
$U_{r,s}({\mathfrak{sl}}_n)$,  such algebras also have an action
of $U_{r,s}({\mathfrak{sl}}_n)$,
and in some cases they admit formal deformations arising from
this action.

In order to obtain these examples, we found it necessary to generalize
the definitions of twisting elements and universal deformation formulas
given in \cite{giaquinto-zhang98}, from elements in bialgebras
to operators on particular categories of modules.
This we do in Section 2.    There are known connections between
twisting elements and $R$-matrices of quasitriangular Hopf algebras
(see Section 2 for the details).
In Section 3,  we factor the $R$-matrix of the quasitriangular
Hopf algebra ${\mathfrak{u}}_{r,s}({\mathfrak{sl}}_{n})$, where $r$ and $s$
are roots of unity (under a mild numerical constraint),
and obtain a twisting element from one of the factors.  
In the special case of  ${\mathfrak{u}}_{1,-1}({\mathfrak{sl}}_2)$,  
we show that this twisting element leads to a universal deformation formula,
and via this alternate approach, essentially recover the example of
\cite{caldararu-giaquinto-witherspoon}.

In Section 4,  we continue under the assumption that $r$ and $s$ are roots
of unity, but
consider the infinite-dimensional quantum group $U_{r,s}({\mathfrak
{sl}}_{n})$. We show that an analogue of an $R$-matrix for $U_{r,s}({\mathfrak
{sl}}_{n})$
 gives rise to a universal deformation formula for certain types of algebras.
In such formulas, it appears to be necessary to work
 with the quantum group $U_{r,s}({\mathfrak
{sl}}_{n})$, or more particularly with $U_{q,q^{-1}}({\mathfrak
{sl}}_{n})$,  rather than with  its  one-parameter
 quotient $U_q({\mathfrak{sl}}_n)$,  
   whose relations do not permit the
types of actions we describe in Section 5.

For the definitions and general theory in Section 2, we assume only that
our algebras are defined over a commutative ring $\kc$.
For the quantum groups in Sections 3 through 5, we assume $\kc$ is a field
of characteristic 0 containing appropriate roots of unity.

\medskip

\section{Twisting elements and deformation formulas}

 Let $B$ be an associative  bialgebra over a  commutative ring $\kc$ with coproduct $\Delta: B \rightarrow B \ot B$ and counit $\epsilon: B \rightarrow \kc$.
An associative  $\kc$-algebra $A$ is called a {\em left $B$-module algebra} if it is a left
$B$-module such that
\begin{equation}\label{modalg1}
  b.1 = \epsilon(b)1, \quad \mbox{and}
\end{equation}
\begin{equation}\label{modalg2}
  b.(aa') = \sum_{(b)}  (b_{(1)}.a)(b_{(2)}.a')
\end{equation}
for all $b\in B$ and $a,a'\in A$.  Here we have adopted the Heyneman-Sweedler convention for
the coproduct,
$\Delta(b) = \sum_{(b)}   b_{(1)} \otimes b_{(2)}$.

Let $\mathcal C$ be a category of $B$-module algebras.
An element $F\in B\otimes B$ is a {\em twisting element for $\mathcal C$}
(based on $B$) if
\begin{equation}\label{twist1}
  (\epsilon\otimes \id)(F) = \id \otimes \id = (\id \otimes \epsilon)(F), \quad  \mbox{and}
\end{equation}
\begin{equation}\label{twist2}
\bigl [(\Delta\otimes \id)(F)\bigr] (F\otimes \id)=\bigl[(\id \otimes\Delta)(F)\bigr]
  (\id \otimes F)
\end{equation}
as operators on the $\K$-modules $A\otimes A$ and $A\otimes A\otimes A$, respectively,
for all objects $A$ in $\mathcal C$, where $\id$ denotes the identity operator.
More generally, we allow $F$ to be a formal infinite sum of elements of
$B\otimes B$, provided that $F$ and the operators on each side of equations (\ref{twist1}) and
(\ref{twist2})  are well-defined operators on the appropriate objects.

The following theorem is essentially the same as
\cite[Thm.~1.3]{giaquinto-zhang98}, but  stated  in terms of the category
$\mathcal C$.    We include a proof for completeness. \medbreak

\begin{thm}\label{af}
Let $A$ be an associative algebra in a category $\mathcal C$ of
$B$-module algebras having multiplication map
$\mu=\mu_A$, and assume  $F$ is a twisting element for $\mathcal C$.
Let $A_F$ denote the $\kc$-module $A$ with multiplication map $\mu\circ F$.
Then $A_F$ is an associative algebra with multiplicative identity
$1=1_A$.
\end{thm}

\begin{proof}
Equations (\ref{modalg1}) and (\ref{twist1}) imply that $1_A$ is
the multiplicative identity in $A_F$. Associativity in $A_F$  is equivalent to
the identity
$$\mu\circ F\circ (\mu\otimes \id)\circ (F\otimes \id) =
  \mu\circ F\circ (\id \otimes\mu)\circ (\id \otimes F)$$
of functions from $A\otimes A\otimes A$ to $A$.
We prove this identity by applying the associativity of $\mu$,
(\ref{twist2}), and (\ref{modalg2}) twice:
\begin{eqnarray*}
\mu\circ F\circ (\mu\otimes \id )\circ (F\otimes \id) & = &
  \mu\circ (\mu\otimes \id)\circ \big [(\Delta\otimes \id)(F)\bigr]\circ (F\otimes \id)\\
    &=& \mu\circ (\id \otimes \mu)\circ \bigl[(\id \otimes \Delta)(F)\bigr]\circ (\id \otimes F)\\
   &=& \mu\circ F\circ (\id \otimes \mu)\circ (\id \otimes F).
\end{eqnarray*}
(This identity may be viewed as a commutative diagram as in
\cite[Thm.~1.3]{giaquinto-zhang98}.)
\end{proof}

Known examples of twisting elements include those arising from
quasitriangular Hopf algebras, and we recall these ideas next,
with a few details for clarity.
A Hopf algebra $H$ is {\em quasitriangular} if the antipode $S$ is bijective,  and
there is an invertible element $R\in H\otimes H$ such that
\begin{equation}\label{tau}
  \tau\bigl(\Delta(h)\bigr)=R\Delta(h)R^{-1}
\end{equation}
for all $h\in H$, where $\tau(a\otimes b)=b\otimes a$,
\begin{eqnarray}
 (\Delta\otimes \id)(R)&=& R^{13}R^{23},\quad \mbox{and}  \label{qt1}\\
(\id \otimes\Delta)(R)&=& R^{13}R^{12}. \label{qt2}
\end{eqnarray}
The notation is standard;   for example,  if $R=\sum_iR^1_i\otimes R^2_i$,
then $R^{13}=\sum_iR^1_i\otimes 1 \otimes R^2_i$.

It can be shown that
$$R^{12}R^{13}R^{23}=R^{23}R^{13}R^{12}$$
(see for example \cite[Prop.~10.1.8]{montgomery93}). Thus left
multiplication of (\ref{qt1}) by $R^{12}=R\otimes 1$ and
left multiplication of (\ref{qt2}) by $R^{23}=1\otimes R$ yields
$$(R\otimes 1)\bigl[(\Delta\otimes \id)(R)\bigr] = (1 \otimes R)\bigl [(\id \otimes\Delta)(R)\bigr],$$
which is very similar to (\ref{twist2}), but applies to {\em right}
$H$-module algebras.
For {\em left} $H$-module algebras,  
$R^{-1}$ (which equals $(S\otimes \id )(R)$ \cite[(10.1.10)]{montgomery93})
satisfies (\ref{twist2}) (as noted in
\cite[pp.~139--40]{giaquinto-zhang98}).    The element
$R_{21} := \tau (R)$  can also
be shown to satisfy (\ref{twist2}) by the following
argument:

Interchange the factors in (\ref{qt1}) and (\ref{qt2})  to obtain
\begin{eqnarray}
(\id \otimes \Delta)(R_{21})&=& R^{21}R^{31},\qquad \mbox{and}  \label{qt3}\\
 (\Delta\otimes \id )(R_{21}) &=&  R^{32}R^{31}.\label{qt4}
\end{eqnarray}

\noindent Using the fact that $R_{21}$ is itself an $R$-matrix for
$H^{\coop}$ (that is, for $H$ with the opposite coproduct), we  
see that
$$R^{32}R^{31}R^{21}=R^{21}R^{31}R^{32}.$$
(This may also be shown directly.)
Thus, multiplying (\ref{qt3}) on the right by $R^{32}$ and multiplying
(\ref{qt4}) on the right by $R^{21}$ yields
\begin{equation}\label{eqn:r21}
\bigl[(\Delta\otimes \id)(R_{21})\bigr](R_{21}\otimes 1) =
\bigl[(\id \otimes \Delta)(R_{21})\bigr](1 \otimes R_{21}).
\end{equation}
Again, as $R_{21}$ is an $R$-matrix  for $H^{\coop}$, it satisfies
(\ref{twist1})
\cite[(10.1.11)]{montgomery93}.
Therefore we have shown the following.\medbreak

\begin{prop}\label{prop:r21}
Let $H$ be a quasitriangular Hopf algebra with \hbox{\rm R}-matrix $R$. Then
$R_{21}$ is a twisting element for left $H$-module algebras.
\end{prop}

Now assume $A$ is an associative $\kc$-algebra with multiplication map $\mu=\mu_A$.
Let $t$ be an indeterminate  and $A[[t]]$  be the algebra
of formal power series with coefficients in $A$.   A map $f\colon A[[t]]
\rightarrow A[[t]]$  is said to have $\deg(f)\geq i$  if 
$f$ takes the ideal $A[[t]]t^j$ to $A[[t]]t^{i+j}$ for each $j$.

A  {\em formal deformation} of $A$
is an associative $\K[[t]]$-algebra structure on the $\kc[[t]]$-module
$A[[t]]$ for which multiplication takes the form
$$\mu_A+t\mu_1+t^2\mu_2+\cdots.$$
Each $\mu_i:A\otimes_{\kc} A\rightarrow A$ is assumed to be a $\kc$-linear map,
extended to be $\K[[t]]$-linear.  Associativity imposes certain
constraints on the maps $\mu_i$;  for example,  $\mu_1$ must be
a {\em Hochschild two-cocycle}:
$$\mu_1(a\ot b)c +\mu_1(ab\ot c) =a\mu_1(b\ot c)+\mu_1(a\ot bc)$$
for all $a,b,c\in A$.  For further discussion of  the conditions on the maps
$\mu_i$, see \cite[p.\ 141]{giaquinto-zhang98}.

Let $B$ be a $\kc$-bialgebra,   ${\mathcal C}$ be a category of
$B[[t]]$-module algebras of the form $A[[t]]$ for $\K$-algebras $A$, and
$F$ be a twisting element for $\mathcal C$.
We say that $F$ is a {\em universal deformation formula for $\mathcal C$}
(based on $B$)   if
\begin{equation} \label{eqn:udf}
   F=1\otimes 1 + \sum_{j\geq 1}F^1_j\otimes F^2_j
\end{equation}
for some elements  $F^1_j,F^2_j$ of $B[[t]]$  such that  
each $F^1_j\otimes F^2_j$  (when  viewed as a
transformation  on $A[[t]]\ot_{\kc[[t]]}A[[t]]$ for any $A[[t]]$ in $\mathcal C$) 
satisfies $\deg (F^1_j\otimes F^2_j)\geq 1$.  We assume 
for each $i \geq 1$ that there
are
only finitely many summands $F_j^1\ot F_j^2$ for which $\deg(F_j^1\ot F_j^2)
\geq i$ but  $\deg(F_j^1\ot F_j^2) \not \geq i+1$.
In this case,  $A[[t]]$ has a (new) associative
algebra structure with multiplication $\mu_{A[[t]]}\circ F$  and
multiplicative identity $1_{A[[t]]} = 1_A$ by Theorem \ref{af}.   In addition, this multiplication
has the following
property:  \ \  for all $a,b\in A$,
$$\mu_{A[[t]]}\circ F(a\otimes b)\equiv \mu_{A}
(a\otimes b) \mod (t).$$
Note that  our definition of a universal deformation formula generalizes that of
\cite[Defn.\ 1.13]{giaquinto-zhang98}, where $F$ is required to have an expression
$1\ot 1 + tF_1+t^2F_2+\cdots$ with each $F_i\in B\ot B$,  and where $\mathcal C$
is the category of $B[[t]]$-module algebras arising from $B$-module
algebras by extension of scalars. The differences are that here $F$ is not required to
satisfy (\ref{twist1}) and (\ref{twist2}) on module algebras outside
of category $\mathcal C$, and
$B$ is not required to take each algebra $A$ in $\mathcal C$ to itself, but
only to $A[[t]]$. We illustrate such actions in Section 5.

The next example is well-known.

\begin{ex}{\rm Let $\K$ be a field of characteristic 0, $B$ be a commutative
$\K$-bialgebra, and $\mathcal C$ be any category of $B[[t]]$-module
algebras of the form $A[[t]]$ where $A$ is a $B$-module algebra. Let $P$ be
the space of primitive elements of $B$ and $p\in P\ot P$. Then
$$
  \exp (tp)=\sum_{i=0}^{\infty}\frac{t^i}{i!} p^i
$$
is a universal deformation formula for $\mathcal C$ (see
\cite[Thm.\ 2.1]{giaquinto-zhang98} for a proof).
}\end{ex}

The following theorem on universal deformation formulas
generalizes a consequence of \cite[Thm.\ 1.3 and Defn.\ 1.13]{giaquinto-zhang98}.
\medbreak

\begin{thm}\label{udf}
Let $B$ be a $\kc$-bialgebra and $\mathcal C$ be a category of
$B[[t]]$-module algebras of the form $A[[t]]$,  where $A$ is a $\K$-algebra.
Let $F$ be a universal deformation formula for $\mathcal C$  based on $B$.
Then for each $A[[t]]$ in $\mathcal C$, $\mu_{A[[t]]}\circ F$ defines a
formal deformation of $A$.
\end{thm}

\begin{proof}
Let $A[[t]]\in {\mathcal C}$.
  By Theorem \ref{af},  $\mu_{A[[t]]}\circ F$ provides an
associative algebra structure on $A[[t]]$.    For $a,a' \in A$,
we may write
$$F^1_j.a=\sum_{k\geq 0}\phi^k_j(a)t^k \ \mbox{ and }
F^2_j.a' =\sum_{l\geq 0}\psi^l_j(a')t^l,$$
where $\phi_j^k:A\rightarrow A$, $\psi_j^l:A\rightarrow A$ are $\kc$-linear
functions.   The
coefficient of $t^i$ in $\mu_{A[[t]]}\circ
(F^1_j\otimes F^2_j)(a\otimes a')$ is
$\sum_{k+l=i}\phi^k_j(a)\psi^l_j(a')$.  As $\deg(F^1_j\otimes F^2_j)\geq 1$,
the coefficient of $t^0$ in $\mu_{A[[t]]}\circ
(F^1_j\otimes F^2_j)(a\otimes a')$
is 0 for each $j\geq 1$.
Therefore we have $\mu_{A[[t]]}\circ F(a\otimes a') =\mu_A(a\otimes a')
+t\mu_1(a\otimes a')+t^2\mu_2(a\otimes a') + \cdots$ where
$$\mu_i(a\otimes a') = \sum_{j\geq 1}\sum_{k+l=i}\phi^k_j(a)\psi^l_j(a').$$
This is necessarily well-defined for each $i$,  since $F$ is assumed to be
a well-defined operator on $A[[t]]\ot_{\K[[t]]}A[[t]]$.
As $F^1_j$, $F^2_j$ act linearly on $A$ for each $j$, the same is true of the
functions $\phi^k_j$, $\psi^l_j$, and thus  $\mu_i$ is bilinear.
\end{proof}  \bigbreak

\section{Twisting elements from finite quantum groups}

Drinfeld doubles of finite-dimensional
Hopf algebras provide a wealth of examples of quasitriangular
Hopf algebras.    In particular, the
two-parameter restricted quantum groups ${\mathfrak{u}}_{r,s}
({\mathfrak{sl}}_{n})$ of \cite{benkart-witherspoon3} are Drinfeld doubles
under some mild assumptions on $r$ and $s$.  We discuss implications
of this
  for twisting and deforming algebras later on.   The restricted quantum group
${\mathfrak{u}}_{r,s}({\mathfrak{sl}}
_{n})$ is a quotient of the unital associative
$\kc$-algebra $U_{r,s}(\mathfrak{sl}_{n})$, which we introduce next.  

Let $\e_1, \dots, \e_n$ denote an orthonormal basis of a Euclidean space
$E={\mathbb R}^n$ with an inner product $\langle \ , \ \rangle$.  
Set $\Pi =
\{\alpha_j = \e_{j}- \e_{j+1} \mid j = 1, \dots, n-1\}$ and
$\Phi =
\{\e_i -\e_j \mid 1 \leq i \neq j \leq n\}$.  Then $\Phi$ is a finite
root system of type A$_{n-1}$ with $\Pi$ a base of simple roots.
Let $\kc$ be a field.
Choose $r,s\in \kc^{\times}$ with $r\neq s$.
The algebra $U=U_{r,s}(\mathfrak{sl}_{n})$  is the unital associative
$\kc$-algebra  
generated by $e_j, \ f_j \  (1 \leq j < n)$,
and $\ \w_i^{\pm 1}, \ (\w_i')^{\pm 1}$ \  ($1 \leq i < n$), subject
to the following relations. \medbreak
 
 \begin{itemize}
\item[{\rm (R1)}]  The $\w_i^{\pm 1}, \ (\w_j')^{\pm 1}$ all commute with one
another and $\w_i \w_i^{-1}= \w_j'(\w_j')^{-1}=1,$\m
\item[{\rm (R2)}] $ \w_i e_j = r^{\la \e_i,\alpha_j\ra }s^{\la
\e_{i+1},\alpha_j\ra}e_j \w_i$ \ \ and \ \ $\w_if_j =
r^{-\la\e_i,\alpha_j\ra}s^{-\la\e_{i+1},\alpha_j\ra} f_j\w_i,$ \m
\item[{\rm (R3)}]
$\w_i'e_j = r^{\la\epsilon_{i+1},\alpha_j\ra}s^{\la\e_{i},\alpha_j\ra}e_j
\w_i'$ \ \ and \ \ $\w_i'f_j =
r^{-\la\epsilon_{i+1},\alpha_j\ra}s^{-\la\e_{i},\alpha_j\ra}f_j \w_i'$, \m
\item[{\rm (R4)}]
$\displaystyle{[e_i,f_j]=\frac{\delta_{i,j}}{r-s}(\w_i-\w_i').}$
\m \item[{\rm (R5)}] $[e_i,e_j]=[f_i,f_j]=0 \ \ \text{ if }\ \
 |i-j|>1, $ \m \item[{\rm (R6)}]
 $e_i^2e_{i+1}-(r+s)e_ie_{i+1}e_i+rse_{i+1}e_i^2 = 0,$
 
 \smallskip
 \hskip -.2 truein $e_i e^2_{i+1} -(r+s)e_{i+1}e_ie_{i+1} +rs
 e^2_{i+1}e_i = 0,$ \m \item[{\rm (R7)}]
 $f_i^2f_{i+1}-(r^{-1}+s^{-1})f_if_{i+1}f_i +r^{-1}s^{-1}f_{i+1}f_i^2 =
 0,$
 
 \smallskip
 \hskip -.2 truein $f_i f^2_{i+1}
 -(r^{-1}+s^{-1})f_{i+1}f_if_{i+1}+r^{-1}s^{-1} f^2_{i+1} f_i=0.$
 \end{itemize}\m
  
The algebra $U$ becomes a Hopf algebra over $\kc$ with $\omega_i$,
$\omega_i'$ group-like and
\begin{eqnarray*}
\Delta(e_i)=e_i\otimes 1+\omega_i\otimes e_i, &\ \ & \Delta(f_i)=1\otimes f_i +
f_i\otimes\omega_i',\\
\epsilon(e_i)=0, &\ \ & \epsilon(f_i)=0,\\
S(e_i)=-\omega_i^{-1}e_i,&\ \ & S(f_i)=-f_i(\omega_i')^{-1}.
\end{eqnarray*} \m

Let $U^0$ be the group algebra generated by all $\omega_i^{\pm 1}$,
$(\omega_i')^{\pm 1}$, and let $U^{+}$ (respectively, $U^{-}$) be the
subalgebra of $U$ generated by all $e_i$ (respectively, $f_i$).

\medbreak

We define
$$\E_{j,j}=e_j \ \ \mbox{ and } \ \ \E_{i,j}=e_i\E_{i-1,j}-
  r^{-1}\E_{i-1,j}e_i \quad  (i>j),$$
$$\F_{j,j}=f_j \ \ \mbox{ and } \ \ \F_{i,j}=f_i\F_{i-1,j}-
s\F_{i-1,j}f_i \quad (i>j).$$
The algebra $U$ has a triangular decomposition $U\cong U^-
\otimes U^0\otimes U^+$, and as discussed in \cite{benkart-kang-lee}, the subalgebras
$U^+$, $U^-$ respectively have
monomial PBW (Poincar\'e-Birkhoff-Witt)  bases

 \begin{equation}\label{pbw1}
\mathcal E: = \{\E_{i_1,j_1}\E_{i_2,j_2}\cdots\E_{i_p,j_p}\mid (i_1,j_1)\leq
(i_2,j_2)\leq\cdots\leq (i_p,j_p)\mbox{ lexicographically}\},
\end{equation}
  
\begin{equation}\label{pbw2}
\mathcal F:= \{\F_{i_1,j_1}\F_{i_2,j_2}\cdots\F_{i_p,j_p}\mid (i_1,j_1)\leq
(i_2,j_2)\leq\cdots\leq (i_p,j_p)\mbox{ lexicographically}\}.
\end{equation}      

The spaces  $U^{\leq 0} = U^{-}\otimes U^0$
and $U^{\geq 0} = U^0 \otimes U^+$ are Hopf subalgebras of $U$.   We view $U^{\leq 0}$ as
a Hopf algebra with the opposite coproduct $\Delta^{\op}$.     Then there is
a Hopf pairing  $(\,|\,): U^{\leq 0}\times U^{\geq 0}  \rightarrow \mathbb K$  given by
\begin{equation}
\begin{array}{ccc}  (f_i\,|\,e_j) &=&
\displaystyle{\frac{\delta_{i,j}}{s-r}}\hspace{.5 in}\\
(\omega_i'\,|\,\omega_j)
&=& r^{\langle \epsilon_j,\alpha_i\rangle}
s^{\langle\epsilon_{j+1},\alpha_i\rangle}\end{array}\label{pairgens} \end{equation}
for $1 \leq i,j < n$, and all other
pairings between generators are 0.         Pairings
between more complicated expressions can be computed from those
values by applying the coproducts $\Delta$ of $U^{\geq 0}$ and $\Delta^{\op}$ of $U^{\leq 0}$
according to the following rules:
 
\begin{equation}\label{pairprops}\end{equation} \vspace{-.35in}\begin{itemize}
\item[{\rm (i)}]  $(1\,|\,b)=\epsilon(b), \ \ \ (b'\,|\,1) =
\epsilon(b')$,   \item[{\rm (ii)}]
$(b'\,|\,bc)=(\Delta^{\op}(b')\,|\,b\otimes c) = \sum_{(b')}  (b_{(2)}'\,|\,b)(b_{(1)}'\,|\,c)$ \\
(where $\Delta(b') = \sum_{(b')} b_{(1)}' \ot b_{(2)}'$), 
\item[{\rm (iii)}]  $(b'c'\,|\,b)=(b'\otimes c'\,|\,\Delta(b))  
= \sum_{(b)}  (b'\,|\,b_{(1)}) (c'\,|\, b_{(2)})$,   \end{itemize}  for all $b,c\in U^{\geq 0}$ and
$b',c'\in U^{\leq 0}$.
\bigbreak

{\em For the rest of this section, we assume that $r$ and $s$ are
roots of unity. Let $r$ be a primitive $d$th root of unity, $s$ be a
primitive $d'$th root of unity, and $\ell$ be the least common multiple
of $d$ and $d'$.  We also assume that $\kc$ contains a primitive
$\ell$th root of unity  $\theta$ and that  $r=\theta^y$, $s=\theta^z$,  where
$y,z$ are nonnegative integers.}
\bigbreak
 
It is shown in \cite{benkart-witherspoon3} that all
$\E^{\ell}_{i,j}$, $\F^{\ell}_{i,j}$, $\omega_i^{\ell}-1$, and
$(\omega'_i)^{\ell}-1$ ($1\leq j \leq i< n$) are central
in $U_{r,s}({\mathfrak{sl}}_{n})$.   The ideal $I_{n}$  
generated by these elements is in fact a Hopf ideal  
\cite[Thm.\ 2.17]{benkart-witherspoon3}.   This is evident from the
following expressions for the coproduct  given in
  \cite[ (2.24) and
   the ensuing text]{benkart-witherspoon3}:
 \begin{equation}\label{eqn:delta-el}
  \Delta(\E^{\ell}_{i,j})=\E^{\ell}_{i,j}\ot 1 +\omega_{i,j}^{\ell}
   \ot\E^{\ell}_{i,j} +s^{\ell (\ell -1)/2}(1-r^{-1}s)^{\ell}\sum
   _{p=j}^{i-1}\E^{\ell}_{i,p+1}\omega_{p,j}^{\ell}\ot\E^{\ell}_{p,j}\hspace{.475in}
\end{equation}
\begin{equation}\label{eqn:delta-fl}
  \Delta(\F^{\ell}_{i,j})= 1\ot \F^{\ell}_{i,j} +\F^{\ell}_{i,j}\ot
   (\omega'_{i,j})^{\ell} +r^{-\ell (\ell-1)/2}(1-r^{-1}s)^{\ell} \sum
    _{p=j}^{i-1}\F_{p,j}^{\ell}\ot \F_{i,p+1}^{\ell}(\omega'_{p,j})^{\ell}
\end{equation}
where $\omega_{p,j}=\omega_p\omega_{p-1}\cdots \omega_j$ and
$\omega_{p,j}'=\omega_p'\omega_{p-1}'\cdots \omega_j'$.
As $I_n$ is a Hopf ideal, the quotient

\begin{equation}\label{resqg}{\mathfrak{u}}:={\mathfrak{u}}_{r,s}({\mathfrak{sl}}_{n}):=
  U_{r,s}({\mathfrak{sl}}_{n})/I_{n}\end{equation}
is a Hopf algebra, called the {\em restricted two-parameter quantum group}.
Furthermore, $\mathfrak{u}$ is finite-dimensional,  as can be readily seen from the PBW-bases
(\ref{pbw1}) and (\ref{pbw2}).   The algebra $U$ is graded by
the root lattice $Q$ of  ${\mathfrak{sl}}_{n}$ by assigning  
$$
  \deg(e_i) = \alpha_i, \ \deg(f_i)=-\alpha_i, \
    \mbox{ and } \ \deg(\omega_i) = 0 = \deg(\omega_i').
$$
Since the generators of $I_n$ are homogeneous, $\mathfrak u$ inherits the
grading.  Note that by (\ref{pairgens}) and (\ref{pairprops}), if $X\in \E$
and $Y\in \F$ with $(Y\, | \, X)\neq 0$, then $\deg(XY)=0$.
 \medbreak

Let $\mathcal E_\ell$ denote the set of monomials in $\mathcal E$ having each $\mathcal E_{i,j}$
appear with exponent at most $\ell-1$.   Identifying  cosets in $\mathfrak u$ with their
representatives, we may assume $\mathcal E_{\ell}$ is a basis for the subalgebra of
$\mathfrak u$ generated by the elements $e_i$.

The following proposition will allow us to use the pairing (\ref{pairgens})
on the quotient algebra $\mathfrak u$.
\medbreak

\begin{prop} \label{radical}
The ideal $I_n$ is contained in the radical of the pairing $(\, |\, )$ of
$U_{r,s}({\mathfrak{sl}}_n)$.  Thus there is an induced pairing on
the quotient ${\mathfrak{u}}_{r,s}({\mathfrak{sl}}_n)$.
\end{prop}

\begin{proof}
For the group-like elements, note that for each pair $i,j$,
\begin{eqnarray*}
  (\omega_i'\, | \,\omega_j^{\ell}-1) &=& (\omega_i'\, |\,\omega_j^{\ell})-
              (\omega_i'|1)\\
           &=& (r^{\langle \epsilon_j,\alpha_i\rangle} s^{\langle
       \epsilon_{j+1},\alpha_i\rangle})^{\ell} -1 = 0,
\end{eqnarray*}
as $r$ and $s$ are $\ell$th roots of 1. Thus $\omega_j^{\ell}-1$
is in the radical of $(\, |\, )$, and similarly for $(\omega_i')^{\ell}-1$.

Now consider $\E_{i,j}^{\ell}$ and let
$X$ be any monomial in $f_1,\ldots, f_{n-1}$ of degree
$-\ell\alpha_i-\cdots -\ell\alpha_j$. Such a monomial is the only type that
potentially has a nonzero pairing with $\E^{\ell}_{i,j}$. As $i\geq j$
and $\ell\geq 2$, we may write $X=Yf_k$ for
some $k$. By (\ref{pairprops})(iii) and (\ref{eqn:delta-el}),
$$
  (X\,|\,\E^{\ell}_{i,j})= (Y\ot f_k \,|\, \E^{\ell}_{i,j}\ot 1 + \omega_{i,j}^{\ell}
    \ot\E^{\ell}_{i,j} +s^{\ell (\ell-1)/2}(1-r^{-1}s)^{\ell}\sum_{p=j}
    ^{i-1}\E^{\ell}_{i,p+1}\omega_{p,j}^{\ell}\ot\E^{\ell}_{p,j}).
$$
Again as $\ell\geq 2$ and $\deg f_k=-\alpha_k$, each term above is 0.
\end{proof}

We will need to apply  the following results from \cite{benkart-witherspoon3},
where the corresponding pairing differs from (\ref{pairgens}) by nonzero
scalar multiples. This difference does not affect the results. (See the
proof of \cite[Thm.\ 4.8]{benkart-witherspoon3} where a relevant
adjustment is made.)
\medbreak

\begin{prop}\label{prop:pair} {\rm \cite[(5.8) and Lem.~4.1]{benkart-witherspoon3}}  
Let $\theta$ be a primitive $\ell$th root
of unity,  and suppose $r = \theta^y$ and $s = \theta^z$.  
Let $\mathfrak b$ be the subalgebra of $\mathfrak u$
generated by $\omega_i, e_i$ $(1 \leq i < n)$ and $\mathfrak b'$ be the subalgebra
of $\mathfrak u$ generated by $\omega_i', f_i$ $(1 \leq i < n)$.   Then the Hopf pairing
$(\,|\,)$ on ${\mathfrak{b}}'\times {\mathfrak{b}}$
satisfying (\ref{pairgens}) and (\ref{pairprops}) is
nondegenerate if
\begin{equation}\label{relprime} (y^{n-1}-y^{n-2}z+-\cdots +(-1)^{n-1}z^{n-1},\ell)=1,
\end{equation}  
(that is, the first expression in the parentheses is relatively prime to $\ell$, the
least common multiple of the orders of $r$ and $s$ as roots of 1).  \end{prop}
\medbreak

We will use the definition of the Drinfeld double $D({\mathfrak b})$ given in \cite{benkart-witherspoon3}:
$D({\mathfrak b})={\mathfrak b}\ot ({\mathfrak b}^*)^{\coop}$ as a coalgebra,
where $({\mathfrak b}^*)^{\coop}$ denotes the dual Hopf algebra with the opposite
coproduct. As an algebra, ${\mathfrak b}$ and $({\mathfrak b}^*)^{\coop}$ become
subalgebras of $D({\mathfrak b})$ under identifications with ${\mathfrak b}\ot 1$
and $1\ot {\mathfrak b}^*$, respectively, and
$$(1\ot b)(a\ot 1)=\sum_{(a),(b)} (b_ {(1)} \, | \, S^{-1}(a_{(1)})) (b_ {(3)}
\, | \, a_ {(3)}) a_ {(2)}\ot b_ {(2)}.$$

\begin{prop}\label{double}
{\rm \cite[Thm.\ 4.8]{benkart-witherspoon3}} Assume  $r = \theta^y$ and $s = \theta^z$,  
where $\theta$ is a primitive $\ell$th root of unity,  and suppose that (\ref {relprime}) holds. Then there is an
isomorphism of Hopf algebras ${\mathfrak{u}}_{r,s}({\mathfrak{sl}}_{n}) \cong D({\mathfrak{b}})$,
where $D({\mathfrak{b}})$ is the Drinfeld double of the Hopf subalgebra ${\mathfrak{b}}$ of
${\mathfrak{u}}_{r,s}(\mathfrak{sl}_{n})$ generated by $\omega_i,e_i \ (1\leq i < n)$.  \end{prop} \m
 
In general, whenever  ${\mathfrak {u}}\cong {\mathfrak{u}}_{r,s}(
{\mathfrak{sl}}_n)\cong D({\mathfrak{b}})$, then $\mathfrak u$ is
quasitriangular with $R$-matrix,
$$
  R=\sum b\ot b^*,
$$
where $b$ runs over a basis of $\mathfrak b$ and $b^*$ runs over the
dual basis of the dual space ${\mathfrak {b}}^*$, which can be identified
with the Hopf subalgebra ${\mathfrak{b}}'$ of $\mathfrak u$ generated
by $\omega_i', f_i$ ($1\leq i <n$) but with the opposite coproduct
$\Delta^{\op}$ \cite[Lem.\ 4.1]{benkart-witherspoon3}.
 
To illustrate this result in a very special case, take  $r = q$, a primitive $\ell$th root of
unity, and $s = q^{-1}$.  Then  $\mathfrak u= {\mathfrak{u}}_{q,q^{-1}}({\mathfrak{sl}}_{n})$ is isomorphic to the  Drinfeld double of
$\mathfrak b$ when $n$ and $\ell$ are relatively prime.   The quotient of
$\mathfrak u$ by the ideal generated by the elements $\omega_i'-\omega_i^{-1}$ for $1 \leq i < n$
is a  one-parameter restricted quantum group related to those which have played
a significant role in the study of algebraic groups in the case
 $\ell = p$, a prime (see for example,  \cite{andersen-jantzen-soergel}).  

\begin{lem}\label{lem:factor}  Assume ${\mathfrak{u}}_{r,s}({\mathfrak{sl}}_{n})$ is
a Drinfeld double as in Proposition \ref{double}.  Then its $R$-matrix factors,
 \begin{equation} R=R_{e,f}R_{\omega,\omega'},  \end{equation}
with $R_{e,f}=\sum \varepsilon \otimes \varepsilon^*$, where $\varepsilon$ runs over the basis
$\mathcal E_\ell$ of $\ell$-power truncated monomials, and $R_{\omega,\omega'}=\sum  w \otimes w^*$,
where $w$ runs over the basis $\Omega: = \{\omega^{\underline c} = \omega_1^{c_1} \cdots
\omega_{n-1}^{c_{n-1}} \mid 0 \leq c_i < \ell$ for all $i\}$  of the group algebra generated
by the $\omega_i \ (1\leq i <  n)$.
\end{lem}

\begin{proof}   Let $\{\varepsilon^* \mid \varepsilon  \in \mathcal E_\ell\}$ denote elements
of ${\mathfrak b}'$ dual to  the vectors of $\mathcal E_\ell$
with respect to this pairing so that $((\varepsilon')^*\mid \varepsilon) = \delta_{\varepsilon',
\varepsilon}$ for all $\varepsilon',\varepsilon \in \mathcal E_\ell$.
For example, if $\varepsilon= \mathcal E_{1,1}\mathcal E_{2,2} = e_1 e_2$, then
$\varepsilon^*$ is a linear combination of $\mathcal F_{1,1}\mathcal F_{2,2} =f_1 f_2$
and $\mathcal F_{2,1} =  f_2f_1 - sf_1 f_2$.      Similarly,
let $\{w^* \mid w \in \Omega\}$ denote the basis of $\kc G'$ dual to $\Omega$,
where $G'$ is the group generated by the $\omega_i'$.    By the triangular
decomposition of $\mathfrak u_{r,s}({\mathfrak {sl}}_n)$,  
the elements $\{ \varepsilon w \mid \varepsilon \in \mathcal E_\ell,
w \in \Omega\}$ form a basis of $\mathfrak b$ (see for example, \cite[(2.16)]{benkart-witherspoon3}).  
Moreover using  $\Delta(e_i) = e_i \ot 1 + \omega_i \ot e_i$
and $\Delta^{\op}(f_i) = f_i \ot 1 + \omega_i' \ot f_i$, the properties of
the pairing in (\ref{pairprops}),  and degree considerations,  we have
\begin{eqnarray*} ((\varepsilon')^*(w')^* \mid \varepsilon w) &=& ((\varepsilon')^*
\ot (w')^* \mid \Delta(\varepsilon w) ) \\
&=& ((\varepsilon')^* \ot (w')^* \mid \Delta(\varepsilon) (w  \ot w) ) \\
&=& ((\varepsilon')^* \mid \varepsilon w)((w')^* \mid w) \\
&=& (\Delta^{\op}((\varepsilon')^*) \mid \varepsilon\ot w) \delta_{w',w} \\
&=& ((\varepsilon')^* \mid \varepsilon) (1 \mid w)  \delta_{w',w} \\
&=& \delta_{\varepsilon',\varepsilon} \delta_{w',w} .
\end{eqnarray*}
This calculation shows that $\varepsilon^* w^*$ is the dual basis element of $\varepsilon w$ relative
to the pairing.  Thus,  
  \begin{eqnarray*}  R  &=&  \sum_{\varepsilon \in \mathcal E_\ell, w \in \Omega}  \varepsilon w \ot  
\varepsilon^*w^* \\
&=& \left(\sum_{\varepsilon \in \mathcal E_\ell} \varepsilon \ot \varepsilon^* \right)\left(
\sum_{w \in \Omega}w \ot w^*\right).    \end{eqnarray*}   \end{proof}

It follows from Lemma \ref{lem:factor} that  $R_{21}=R_{f,e}
R_{\omega',\omega}$, where $R_{f,e} = \sum_{\varepsilon \in \mathcal E_\ell} \varepsilon^* \ot \varepsilon$ = $(R_{e,f})_{21}$
and $R_{\omega',\omega} = \sum_{w \in \Omega}  w^* \ot w = (R_{\omega,\omega'})_{21}$.
Equation (\ref{eqn:r21}) holds for $R_{21}$, that is
\begin{equation}\label{ref}
\bigl[(\Delta\otimes \id )(R_{f,e}R_{\omega',\omega})\bigr](R_{f,e}R_{\omega',
\omega}\otimes 1 ) = \bigl[(\id \otimes\Delta)(R_{f,e}R_{\omega',\omega})\bigr]
(1 \otimes R_{f,e}R_{\omega',\omega}).
\end{equation}
We will show that the factor $R_{f,e}$ alone satisfies this equation and so is
a twisting element.

Let $\mathcal{C}$ be the category of left $\mathfrak{u}$-module algebras.  
Let $G$ be the group generated by the $\omega_i$ and $G'$ the group generated
by the $\omega_i'$ ($1\leq i<n$).
Since  $\omega_i,\omega_i' \ (1 \leq i < n)$  generate the finite abelian group
$G\times G'$ of order a power of $\ell$,  and $\K$ contains a primitive $\ell$th root of 1,
each $A\in {\mathcal{C}}$ decomposes into a direct sum of
common eigenspaces (or weight spaces) for the $\omega_i,\omega_i'$.
That is, $A=\oplus_{\chi}A_{\chi}$ for some group homomorphisms
$\chi:G\times G'\rightarrow \K^{\times}$, where $\omega_i . a =\chi(\omega_i)a$
and $\omega_i'.a =\chi(\omega_i')a$ for all $a\in A_{\chi}$, $i=1,\cdots, n-1$.
We will use this decomposition to show that
$R_{f,e}$ is itself a twisting element for ${\mathcal {C}}$. \medbreak

\begin{thm}\label{thm:twist}  Let ${\mathfrak{u}}:={\mathfrak{u}}_{r,s}({\mathfrak{sl}}_{n})$.  
The factor $R_{f,e}$ of the opposite $R_{21}$ of the $R$-matrix for
$\mathfrak u$ is a twisting element for the category of left $\mathfrak u$-module
algebras.
\end{thm}    

\begin{proof}
Let $A$ be a left $\mathfrak u$-module algebra, and apply both sides of (\ref{ref})
to $a\otimes a'\otimes a''\in A_{\chi}\otimes A_{\psi}\otimes
A_{\phi}$, for weights $\chi,\psi,\phi$ of $A$.
Note that $R_{\omega',\omega}$ acts as the  scalar  $\xi(\chi,\psi) := \sum_{w \in \Omega}
\chi(w^*)\psi(w)$ on $a\otimes a'$.    Applying the left
side of (\ref{ref}) to $a\otimes a'\otimes a''$, we have
\begin{eqnarray*}&&\bigl[(\Delta\otimes \id )(R_{f,e}R_{\omega',\omega})\bigr](R_{f,e}\otimes 1)
\xi(\chi,\psi) a\otimes a'\otimes a'' \\
&& \hspace{1.2truein}  = \bigl[(\Delta\otimes \id )(R_{f,e})\bigr](R_{f,e}\otimes 1)
  \xi(\chi\cdot\psi,\phi)\xi(\chi,\psi)a\otimes a'\otimes a'', \end{eqnarray*}
since each term of $R_{f,e}$ when applied to $a\otimes a'$
changes the weight of $a$ and the weight of $a'$ by multiplication by functions
which are inverses of one another, with a net effect of no change at all in the weight.  Similarly,
applying the right side of (\ref{ref}) to $a\otimes a'\otimes a''$,
we have
\begin{eqnarray*}&&\bigl[(\id \otimes\Delta)(R_{f,e}R_{\omega',\omega})\bigr]
(1\otimes R_{f,e})
\xi(\psi,\phi)a\otimes a'\otimes a'' \\
&& \hspace{1.2truein}  = \bigl[(\id \otimes\Delta)(R_{f,e})] (1\otimes R_{f,e})
\xi(\chi,\psi\cdot\phi)\xi(\psi,\phi) a\otimes a'\otimes a''.\end{eqnarray*}
Now $R_{\omega',\omega}$ is itself an $R$-matrix for the group algebra $\kc [G\times G']$,
because our constructions show that this group algebra is the
Drinfeld double of $\kc[G]$.    Similar calculations to those
above prove that
$$\xi(\chi\cdot\psi,\phi)\xi (\chi,\psi)= \xi(\chi,\psi\cdot\phi)
\xi(\psi,\phi)$$
for all weights $\chi,\psi,\phi$. (Note this implies
that $\xi$ is a 2-cocycle for the dual group to $G\times G'$.)
Moreover, as $R_{\omega',\omega}$ is invertible, none of these
values is zero. Thus we may cancel the factors
$\xi(\chi\cdot\psi,\phi)\xi(\chi,\psi)$ and $\xi(\chi,\psi\cdot\phi)
\xi(\psi,\phi)$ from the above expressions and apply (\ref{ref})
to obtain
$$[(\Delta\otimes \id)(R_{f,e})](R_{f,e}\otimes 1)
(a\otimes a'\otimes a'')=[(\id\otimes\Delta)(R_{f,e})]
(1\otimes R_{f,e})(a\otimes a'\otimes a''),$$
as desired. We also have $(\epsilon\otimes 1)(R_{f,e})
=1\otimes 1=(1\otimes\epsilon)(R_{f,e})$ by the
definition of $R_{f,e}$.
\end{proof}

\medbreak
\begin{rem}{\rm Theorem \ref{deltaf} in the next section provides an alternate
approach to proving Theorem \ref{thm:twist}.}  \end{rem}
\medbreak

\begin{exs}{\rm We will give some examples of $\mathfrak u$-module algebras, which 
may be twisted by $R_{f,e}$ according to Theorems \ref{af} and \ref{thm:twist}.
First let $V$ be the natural $n$-dimensional module for $U=U_{r,s}
({\mathfrak{sl}}_n)$;  that is,  $V$ has a basis $v_1,\ldots,v_n$ and
\begin{eqnarray*}
  e_i . v_j &=& \delta_{i,j-1}v_{j-1}\\
  f_i . v_j &=& \delta_{i,j}v_{j+1}\\
  \omega_i . v_j &=& r^{\delta_{i,j}} s^{\delta_{i,j-1}}v_j,   \\
  \omega_i' . v_j &=& r^{\delta_{i,j-1}}s^{\delta_{i,j}}v_j,   \\
\end{eqnarray*}
where $v_0 = 0 =v_{n+1}$.   Because  each ${\mathcal {E}}^{2}_{i,j}$ acts as 0,
 $V$ becomes a
$\mathfrak u$-module via the induced action.
 This action extends to a representation of $\mathfrak u$ on
the tensor algebra $T(V)=\bigoplus_{k=0}^{\infty}V^{\ot k}$, where
the action on $V^{\ot k}$ is by $\Delta^{k-1}$ for $k\geq 1$.
When $k=0$, then $V^{\ot k}=\kc$, and the action is given by the
counit. By definition then, $T(V)$ is a $\mathfrak u$-module algebra.

Let $J$ be the (two-sided) ideal of $T(V)$ generated by all elements
of the form
$$v_i\ot v_j - rv_j\ot v_i \ \ \ (j>i).$$
Then $J$ is homogeneous, $J=\bigoplus_{k=2}^{\infty}J_k$, where
$J_k=J\cap V^{\ot k}$. In fact,  $J$ is a $\mathfrak u$-submodule of $T(V)$:
A computation such as in \cite[Prop.\ 5.3]{benkart-witherspoon2}
shows that $J_2:=J\cap V^{\ot 2}$ is a $\mathfrak u$-submodule
of $V^{\ot 2}$, and as $J=T(V)\ot J_2\ot T(V)$, the result follows.
Thus,  the quantum plane  $\K_r[x_1,\ldots,x_n]$,
 with $x_ix_j = r x_j x_i$ if $j > i$, is isomorphic
 to the $\mathfrak u$-module algebra  $T(V)/J$, under the identification
$x_i=v_i+J$.   Let
$$
  x(\underline{d})=x(d_1,\ldots,d_n)=x_1^{d_1}\cdots x_n^{d_n}.
$$
Then the $\mathfrak u$-module action is given by
\begin{eqnarray*}
  e_i . x(\underline{d}) &=& r^{d_i-d_{i+1}+1}[d_{i+1}]
          x(d_1,\ldots,d_i+1,d_{i+1}-1,\ldots,d_n),\\
  f_i.x(\underline{d})&=& r^{d_{i+1}-d_i +1}[d_i]x(d_1,\ldots
     d_i-1,d_{i+1}+1,\ldots,d_n),\\
  \omega_i.x(\underline{d})&=& r^{d_i}s^{d_{i+1}}x(\underline{d}), \\
  \omega_i'.x(\underline{d})&=& r^{d_{i+1}}s^{d_i}x(\underline{d}),
\end{eqnarray*}
where
$$
  [d]:=\frac{r^d-s^d}{r-s}.
$$

Other quotients of $T(V)$ provide examples as well:    Choose a positive
integer $p$, and let $J$ be the ideal of $T(V)$ generated by
the $\mathfrak u$-submodule of $T(V)$ having as generators all $v_i^{\ot p}$.
As an abbreviation, we will write $v_i^p=v_i^{\ot p}$, and similarly
leave out the tensor product symbol in the notation for words in
$v_1,\ldots,v_n$.
For example, if $n=2$ and $p=2$, the $\mathfrak u$-submodule
of $V^{\ot 2}$ generated by $v_1^{2}$ and $v_2^{2}$ has
basis $v_1^{2}, \ v_2^{2}, \ v_1v_2+sv_2v_1$.
If $n=2$ and $p=3$, the $\mathfrak u$-submodule of $V^{\ot 3}$
generated by $v_1^{3}$ and $v_2^{3}$ has basis $v_1^{3},
\ v_2^{3}, \ v_1v_2^2 +sv_2v_1v_2+s^2v_2^2v_1, \
v_1^2v_2+sv_1v_2v_1+s^2v_2v_1^2$.
If $q$ is a primitive third root of 1, and $r=q$, $s=q^{-1}$, then
$\mathfrak u = {\mathfrak{u}}_{q,q^{-1}}({\mathfrak{sl}}_2)$ acts on the quotient of $T(V)$ by the ideal  
generated by $v_1^2v_2+sv_1v_2v_1 + s^2v_2v_1^2$,  $v_1v_2^2 + sv_2 v_1v_2 + s^2v_2^2v_1$,
which is the down-up algebra $A(1+s^2, -s^2,0) = A(1+r,-r,0)$
in the notation of \cite{benkart-roby}.  
 This example is related to one appearing in the work of Montgomery
and Schneider \cite{montgomery-schneider01} on skew derivations and actions of the double of
the Taft algebra.

The family of down-up algebras $A(1+r,-r,0)$  (where
$r$ is an arbitrary  primitive $\ell$th root of unity)  is especially interesting,  as
the finite-dimensional modules for these algebras are completely reducible \cite{carvalho-musson}.  This family
 includes the universal enveloping algebra
$U(\mathfrak {sl}_2)$ for $r =1$, the universal enveloping algebra $U(\mathfrak {osp}_{1,2})$ of
the Lie superalgebra $\mathfrak {osp}_{1,2}$
for $r = -1$, and the algebras appearing in this example for $r$ a primitive third root of 1.}\end{exs}

\medbreak
\begin{ex}  {\rm Here we generalize the last example by allowing $n$ to
be arbitrary.  Again we suppose $r = q$, $s = q^{-1}$ for $q$ a primitive
third root of 1,  but take $\mathfrak u = \mathfrak u_{r,s}(\mathfrak{sl}_n)$ for any $n \geq 3$.      Consider the subspace  $Y$
of $V^{\ot 3}$  spanned  by the following elements:  

\begin{eqnarray*} y(i,j,k) &=& v_iv_jv_k + s v_jv_iv_k + sv_iv_kv_j+s^2 v_jv_kv_i+s^2v_kv_iv_j + v_kv_jv_i \\
&& \hspace{2 truein}
\quad (1\leq i<j<k\leq n), \\
y(i,i,k) &=& v_i^2 v_k + s v_iv_kv_i + s^2 v_kv_i^2\qquad  \qquad (1 \leq i < k\leq n), \\
y(i,k,k) &=& v_iv_k^2 + s v_kv_iv_k + s^2 v_k^2 v_i \ \qquad \qquad (1 \leq i < k\leq n),\\
y(i,i,i)&=& 0\hspace{2 truein} \qquad (1 \leq i \leq n). \end{eqnarray*}

\noindent  From the formulas below it is easy to see that $Y$ is a $\mathfrak u$-submodule of  $T(V)$:

\begin{eqnarray*}  
\omega_t.y(i,j,k) &=& r^{\la \epsilon_t, \epsilon_i+\epsilon_j + \epsilon_k\ra}
s^{\la \epsilon_{t+1}, \epsilon_i+\epsilon_j + \epsilon_k\ra}  y({i,j,k}), \\
\omega_t'.y({i,j,k})&=& r^{\la \epsilon_{t+1}, \epsilon_i+\epsilon_j + \epsilon_k\ra}
s^{\la \epsilon_{t}, \epsilon_i+\epsilon_j + \epsilon_k\ra}  y({i,j,k}),\\
 \end{eqnarray*}
for all $1 \leq i \leq j \leq k \leq n$;   and  using the fact that $r + s = -1$, we have
\begin{eqnarray*}  
e_t.y(t+1,j,k) &=&  y(t,j,k) \ \qquad  \qquad (t+1 \leq j < k)\  \hbox{\rm or} \  (t+1 < j= k), \\
e_t.y(i,t+1,k) & =& (-1)^{\delta_{i,t}} y(i,t,k) \qquad \quad  \qquad (i < t+1 \leq k), \\
e_t.y(i,j,t+1) & = & (-1)^{\delta_{j,t}} y(i,j,t)\qquad  \qquad \qquad (i \leq j \leq t), \\
e_t.y(i,j,k) & = & 0\hspace{2 truein}  \hbox{\rm otherwise}; \end{eqnarray*}
\begin{eqnarray*}  
f_t.y(t,j,k) &=& (-1)^{\delta_{j,t+1}}  y(t+1,j,k)  \qquad  \qquad (t < j \leq  k), \\
f_t.y(i,t,k) & =& (-1)^{\delta_{k,t+1}} y(i,t+1,k)\qquad  \qquad (i \leq t < k), \\
f_t.y(i,j,t) & = & y(i,j,t+1) \quad \quad  \qquad (i \leq j < t) \ \hbox{\rm or} \ (i < j = t),  \\
f_t.y(i,j,k) & = & 0\hspace{2 truein}  \hbox{\rm otherwise.} \end{eqnarray*}
The restricted quantum group  $\mathfrak u$ acts on the
quotient algebra $A = T(V)/ \la Y \ra$  obtained by factoring out the ideal generated by $Y$.
In the quotient, each  pair $v_i, v_k$ ($1 \leq i < k \leq n$) generates a down-up algebra $A(1+s^2, -s^2,0) = A(1+r,-r,0)$.} \end{ex}

\bigbreak

\begin{ex} {\rm As a special instance of Theorem \ref{thm:twist}, consider the algebra
${\mathfrak{u}}_{1,-1}({\mathfrak{sl}}_2)$ for which
$$R_{f,e}=1\otimes 1 -2f\otimes e.$$
(Note we omit subscripts when discussing $\mathfrak {sl}_2$,  as there
is only one of each type of generator.)
We know from Theorem \ref{thm:twist} that this is a twisting
element, but this also may be verified directly quite easily.
This quantum group has as a quotient
$${\mathfrak{u}}_{-1}({\mathfrak{sl}}_2):={\mathfrak{u}}
_{1,-1}({\mathfrak{sl}}_2)/(\omega-\omega'),$$
and in ${\mathfrak{u}}_{-1}({\mathfrak{sl}}_2)$, the images of $e$ and $f$ commute by
relation (R4).
Therefore, the defining relations of ${\mathfrak{u}}_{-1}({\mathfrak{sl}}_2)$ are all
homogeneous with respect to powers of $e$. Consequently,
the following is a twisting element based on
${\mathfrak{u}}_{-1}({\mathfrak{sl}}_2)[[t]]$, as
may be checked directly:
$$F=1\otimes 1 -2tf\otimes e.$$
Thus,  $F$ is a universal deformation formula for the category of
${\mathfrak{u}}_{-1}({\mathfrak{sl}}_{2})[[t]]$-module algebras
arising from ${\mathfrak{u}}_{-1}({\mathfrak{sl}}_2)$-module algebras by extension of scalars.
In fact, there is a one-parameter
family of such deformation formulas, given by $F_c
=1\otimes 1 + c tf\otimes e$, for any $c\in\kc$.
The choice $c=1$ yields precisely the universal deformation
formula of \cite[Lem.~6.2]{caldararu-giaquinto-witherspoon}
that was applied to a certain ${\mathfrak{u}}_{-1}({\mathfrak{sl}}
_2)$-module algebra $A$ (given by a crossed product of a
polynomial ring with a finite group) to obtain a formal
deformation of $A$.

However, it is an accident due to the choices $\ell=2$,
$n=2$ for ${\mathfrak{u}}_{1,-1}({\mathfrak{sl}}_2)$ that this
technique produces a universal deformation formula directly, as
relation (R4) of ${\mathfrak{u}}_{r,s}
({\mathfrak{sl}}_{n})$ is {\em not} homogeneous
with respect to the powers of the $e_i$ (nor of the $f_i$), and in general,
there is no reasonable quotient in which it becomes so.
In the next section we will remedy this
situation by returning to the infinite-dimensional
Hopf algebra $U=U_{r,s}({\mathfrak{sl}}_n)$ and choosing
an {\em action} of $U$ that itself incorporates the
indeterminate $t$. Using this infinite-dimensional Hopf algebra
necessitates a more complicated
(but related) construction of a twisting element.} \end{ex}
\bigbreak

\section{Twisting elements from infinite quantum groups}

We consider  the infinite-dimensional quantum group
$U=U_{r,s}({\mathfrak{sl}}_{n})$, defined in the previous section,
where $r$ and $s$ are roots of unity. We adopt the same notation
as before, so that   $\ell$ is the least common multiple
of the orders of $r$ and $s$ as roots of $1$, and $Q=\oplus_{i=1}^{n-1}\Z\alpha_i$
is the root lattice of ${\mathfrak{sl}}_n$. We will also assume (\ref{relprime})
holds, so that we may use Proposition \ref{prop:pair}.
Let $\KK$ be a commutative $\K$-algebra, and extend $U$ to
a $\KK$-algebra $U_{\KK}$ that is free as a
$\KK$-module. (For the applications in the next section,
$\KK$ will be the $\K$-algebra of Laurent polynomials in $t$ with
finitely many negative powers of $t$, and $U_{\KK}$ will be a similar extension.)

As in \cite{benkart-witherspoon1}, we will be interested in weight
modules, this time for $U_{\KK}$. We will consider only $U_{\KK}$-modules
that are free as $\KK$-modules.
For every $\lambda\in Q$, define the corresponding algebra homomorphism
$\hat{\lambda} :U^0\rightarrow\kc$ by
\begin{equation}\label{lambdahat}
\hat{\lambda}(\omega_j)=r^{\langle \e_j,\lambda\rangle}s^{\langle
\e_{j+1},\lambda\rangle} \ \ \mbox{ and } \ \ \hat{\lambda}(\omega_j')
=r^{\langle\e_{j+1},\lambda\rangle}s^{\langle\e_j,\lambda\rangle},
\end{equation}
which may be extended to yield an algebra homomorphism
from $U^0_{\KK}$ to $\KK$. As $r$ and $s$ are roots of unity,
the set $\widehat{Q}=\{\hat{\lambda}\mid \lambda\in Q\}$ is finite.

Let $M$ be a $U_{\KK}$-module, assumed to be free as a $\KK$-module.
If as a $\UK^0$-module, $M$ decomposes into the direct sum of
eigenspaces
$$M_{\chi}=\{m\in M\mid (\omega_i-\chi(\omega_i))m=0=(\omega_i'
-\chi(\omega_i'))m \mbox{ for all } i\}$$
for algebra homomorphisms $\chi:\UK^0\rightarrow \KK$,
we say $\UK^0$ {\em acts semisimply on $M$}.
The homomorphisms $\chi$ such that $M_{\chi}\neq 0$ are called the
{\em weights} of $M$.  Note that
\begin{equation}\label{ejfj}
e_j.M_{\chi}\subseteq M_{\chi\cdot\widehat{\alpha_j}} \ \
\mbox{ and } \ \ f_j.M_{\chi}\subseteq M_{\chi\cdot(\widehat{-\alpha_j})}
\end{equation}
by (R2), (R3), and (\ref{lambdahat}). It follows that
if  $M$ is a simple $U_{\KK}$-module with a nonzero weight
space $M_\chi$, then $\UK^0$ acts semisimply on $M$, since the (necessarily direct) sum of
the weight spaces
is a submodule.      In fact, if  $M=\UK.m$, any cyclic $\UK$-module generated
by a weight vector $m \in M_\chi$, then $\UK^0$ acts semisimply on $M$,  and
all the weights of $M$ are of the form $\chi\cdot
\hat{\zeta}$ ($\zeta\in Q$).  \medskip

Let $\CN$ be the category of unital $\UK$-modules $M$ that are free
left  $\KK$-modules and satisfy the following conditions:
\medbreak
\begin{itemize}
\item[($\CN1$)]  $\UK^0$ acts semisimply on $M$;
 
\item[($\CN2$)]  For each $i,j$ ($1\leq j\leq i< n$),  both $\E_{i,j}^{\ell}$
and $\F_{i,j}^{\ell}$ annihilate $M$.
\end{itemize}

\medbreak

\noindent
Note that $\CN$ is closed under direct sums and quotients, and it follows from
(\ref{eqn:delta-el}) and (\ref{eqn:delta-fl}) that $\CN$ is closed under tensor
products.

Examples of modules $M$ satisfying ($\CN$2) are the ${\mathfrak{u}}_{r,s}
({\mathfrak{sl}}_{n})$-modules (see the previous section) with
scalars extended to $\KK$,
which can be viewed as $\UK$-modules.
However we will also be interested in some modules in category $\CN$
that do not correspond to ${\mathfrak{u}}_{r,s}({\mathfrak{sl}}_{n})$-modules,
where neither $\omega_i^{\ell}$ nor $(\omega_i')^{\ell}$ acts as  the identity.

We remark that ($\CN$2) has been included so that certain operators
will be well-defined. It is not necessarily
true that $e_i,f_i$ act nilpotently on all finite-dimensional modules
as happens in the non-root of unity case \cite[Cor.\ 3.14]{benkart-witherspoon1}.
The argument used there fails,  as $\widehat{Q}$ is finite.

For any two modules $M$ and $M'$ in $\CN$, we will construct
a $\UK$-module homomorphism $R_{M',M}:M'\otimes M\rightarrow M\otimes M'$
by Jantzen's method (see \cite[Ch.~7]{jantzen96}).
The main difference between the two-parameter version of Jantzen's
method (see \cite[Sec.~4]{benkart-witherspoon1}) and what we are about
to do here is that now we
are assuming both $r$ and $s$ are roots of unity, a case excluded from  
consideration in \cite{jantzen96, benkart-witherspoon1}.
As a consequence, we will not know whether $R_{M',M}$ is invertible. However
we will make the necessary adjustments to show
that most of the arguments of \cite[Secs.~4,5]{benkart-witherspoon1},
in particular those needed to obtain a twisting element, apply in this context.

The desired function $R_{M',M}$ will be the composition of three $\KK$-linear
functions $\tau$, $\tilde{f}$, and $F$, where $\tau$ is the
 map that interchanges tensor factors as before, and $\tilde{f}$ and $F$ are as
follows. (Ultimately we will show that $F$ is a twisting element, accounting
for the choice of notation.)

As $M$ satisfies ($\CN$1), by (\ref{ejfj}) there are algebra homomorphisms
$\chi: U^0_{\KK}\rightarrow \KK$ such that $M=\oplus_{\chi}M(\chi)$, where
$M(\chi)$ is a $U_{\KK}$-submodule having weights contained in $\chi\cdot
\widehat{Q}$. (For example, we may take the sum over a set of representatives
$\chi$ of cosets of $\widehat{Q}$ in ${\rm {Alg}}_{\KK}(U^0_{\KK},\KK)$.)
For the purposes of this section, it will suffice to deal with
each summand $M(\chi)$ separately. Assume $M=M(\chi)$ and $M'=M'(\psi)$
are such modules. Let $\lambda,\mu\in Q$ and $m\in M_{\chi\cdot\hat{\lambda}}$,
$m'\in M_{\psi\cdot\hat{\mu}}$. Set
$$
  \tilde{f}(m\ot m')=f_{\chi,\psi}(\lambda,\mu)(m\ot m')
$$
where $f_{\chi,\psi}:Q\times Q\rightarrow \KK^{\times}$ is given by
$$
  f_{\chi,\psi}(\lambda,\mu)=\psi(\omega_{\lambda}^{-1})\chi(\omega_{\mu}')
    (\omega_{\mu}' \, | \, \omega_{\lambda}^{-1})
$$
with $\omega_{\lambda}=\omega_1^{\lambda_1}\cdots \omega_{n-1}^{\lambda_{n-1}}$,
$\omega_{\mu}'=(\omega_1')^{\mu_1}\cdots (\omega_{n-1}')^{\mu_{n-1}}$,
$\lambda=\sum_{i=1}^{n-1}\lambda_i\alpha_i$, $\mu=\sum_{i=1}^{n-1}\mu_i\alpha_i$,
and the Hopf pairing is as in (\ref{pairgens}), (\ref{pairprops}).
Thus $f_{1,1}$ is the function $f$ of \cite[(4.2)]{benkart-witherspoon1} as
restricted to the root lattice $Q$, and it can be shown that $\tilde{f}$
generalizes the function $\tilde{f}$ there by looking at cosets of
the weight lattice modulo $Q$.
In particular, the following hold:
\begin{equation}
\begin{array}{rcl}
  f_{\chi,\psi}(\lambda +\mu,\nu) &=& f_{\chi,\psi}(\lambda,\nu)\psi
             (\omega_{\mu}^{-1})(\omega_{\nu}'\, | \, \omega_{\mu}^{-1})\\
  f_{\chi,\psi}(\lambda,\mu + \nu) &=& f_{\chi,\psi}(\lambda,\mu)\chi(
            \omega_{\nu}')(\omega_{\nu}' \, | \, \omega_{\lambda}^{-1})\\
  (\omega_{\mu}' \, | \, \omega_{\lambda}^{-1})&=& \hat{\mu}
   (\omega_{\lambda}^{-1}) \ = \ \hat{\lambda}(\omega_{\mu}').
\end{array}\label{eqn:fid}
\end{equation}

The definition of $F$ is similar to that of $\Theta$ in
\cite[Sec.\ 4]{benkart-witherspoon1}.
We will construct $F$ as a sum of elements of $U\ot U$, which then may  
also be considered elements of $U_{\KK}\ot_{\KK}U_{\KK}$. The subalgebra
$U^+$ of $U$ generated by $1$ and $e_i$ ($1\leq i<n$) may be decomposed
as
$$
   U^+=\bigoplus_{\zeta\in Q^+}U^+_{\zeta}
$$
where $U^+_{\zeta} = \{x\in U^+ \mid x \mbox{ is homogeneous of degree }\zeta\}$,
and $Q^+=\bigoplus_{i=1}^{n-1} \Z^{\geq 0}\alpha_i$.
For each $\zeta\in Q^+$, let $\overline{U}_{\zeta}^+$ be the linear
span of all PBW basis elements in $U^+_{\zeta}$ (see (\ref{pbw1})) in which
the power of each $\E_{i,j}$ is less than $\ell$. Similarly define
$\overline{U}_{-\zeta}^-$. By Propositions \ref{radical} and \ref{prop:pair},
the spaces $\overline{U}^+_{\zeta}$, $\overline{U}^-_{-\zeta}$ are
nondegenerately paired under the assumption that (\ref{relprime}) holds.
We define
\begin{equation}\label{eqn:F}
   F=\sum_{\zeta\in Q^+} F_{\zeta}
\end{equation}
where $F_{\zeta}=\sum_{k=1}^{\overline{d}_{\zeta}} v^{\zeta}_k\ot
u^{\zeta}_k$, \ $\overline{d}_{\zeta} = \dim_{\K} \overline{U}^+_{\zeta}$,
\ $\{u^{\zeta}_k\}_{k=1}^{\overline{d}_{\zeta}}$ is a basis for
$\overline{U}^+_{\zeta}$,  and $\{v^{\zeta}_k \}_{k=1}^{\overline{d}_{\zeta}}$
the dual basis for $\overline{U}^-_{-\zeta}$. Note that if
$\overline{U}^+_{\zeta} = 0$, then $F_{\zeta}=0$, and if
$\zeta\in Q\setminus Q^+$, we will also set $F_{\zeta}=0$ for convenience.

As $\Delta(e_i)=e_i\ot 1 + \omega_i\ot e_i$, for all $x\in U_{\zeta}^+$
we have
$$
   \Delta(x)\in\sum_{0\leq\nu\leq\zeta} U^+_{\zeta-\nu}\omega_{\nu}
      \ot U^+_{\nu},
$$
where $\nu\leq \zeta$ means $\zeta-\nu\in Q^+$.
For each $i$, there are elements $p_i(x)$ and $p_i'(x)\in U^+_{\zeta-\alpha_i}$
such that
\begin{eqnarray*}
  \Delta(x) &=& x\ot 1 + \sum_{i=1}^{n-1} p_i(x)\omega_i\ot e_i + \pi(x),\\
   \Delta(x) &=& \omega_{\zeta}\ot x + \sum_{i=1}^{n-1} e_i\omega_{\zeta-\alpha_i}
       \ot p_i'(x) + \pi '(x),
\end{eqnarray*}
where $\pi(x)$ (respectively, $\pi '(x)$) is a sum of terms involving products
of more than one $e_j$ in the second factor (respectively, in the first factor).
Similarly, if $y\in U^-_{-\zeta}$, we define $p_i(y)$ and $p_i'(y)$ by
\begin{eqnarray*}
  \Delta(y) &=& y\ot \omega_{\zeta}' + \sum_{i=1}^{n-1}p_i(y)\ot f_i\omega_
    {\zeta -\alpha_i}' + \pi(y),\\
  \Delta(y) &=& 1\ot y + \sum_{i=1}^{n-1} f_i\ot p_i'(y)\omega_i' +\pi'(y).
\end{eqnarray*}
The following identities from \cite{benkart-witherspoon1} hold in our context,
as the proof consists of calculations that do not use properties of the
parameters $r$ and $s$.

\begin{lem} {\rm{\cite[Lem.\ 4.6, Lem.\ 4.8]{benkart-witherspoon1}}}
\label{lem:xy}
For all $x\in U^+_{\zeta}$, $x'\in U^+_{\zeta'}$, and $y\in U^-$, the
following hold:
\begin{itemize}
\item[(i)] $(f_iy \, | \, x)=(f_i \, | \, e_i)(y\, | \, p_i'(x))=
         (s-r)^{-1}(y\, |p \, _i'(x))$.
\item[(ii)] $(yf_i \, | \, x)=(f_i \, | \, e_i)(y,p_i(x))=(s-r)^{-1}(y\, | \,
     p_i(x))$.
\item[(iii)] $f_ix-xf_i= (s-r)^{-1}(p_i(x)\omega_i-\omega_i'p_i'(x))$.
\end{itemize}
For all $y\in U^-_{-\zeta}$, $y'\in U^-_{-\zeta'}$, and $x\in U^+$, the
following hold:
\begin{itemize}
\item[(iv)] $(y\, | \, e_ix) = (f_i \, | \, e_i)(p_i(y)\, | \, x)
     =(s-r)^{-1}(p_i(y)\, | \,x)$.
\item[(v)] $(y,xe_i) = (f_i\, | \, e_i)(p_i'(y)\, | \, x)=(s-r)^{-1}
            (p_i'(y)\, | \, x)$.
\item[(vi)] $e_iy-ye_i= (s-r)^{-1}(p_i'(y)\omega_i'-\omega_ip_i(y))$.
\end{itemize}
\end{lem}

Since  the spaces $\overline{U}^+_{\zeta}$ and $\overline{U}^-_{-\zeta}$
are nondegenerately paired with dual bases $\{u^{\zeta}_k\}_{k=1}^{\overline{d}_{\zeta}}$
and $\{v^{\zeta}_k\}_{k=1}^{\overline{d}_{\zeta}}$, for each $x\in
\overline{U}^+_{\zeta}$ and $y\in \overline{U}^-_{-\zeta}$, we have
\begin{equation}\label{eqn:lincomb}
  x=\sum_{k=1}^{\overline{d}_{\zeta}} (v^{\zeta}_k\, | \, x) u^{\zeta}_k
   \ \mbox{ and } \ y=\sum_{k=1}^{\overline{d}_{\zeta}}(y\, | \, u^{\zeta}_k)
    v^{\zeta}_k.
\end{equation}
The next lemma is a modified version of \cite[Lem.\ 4.10]{benkart-witherspoon1} or
\cite[Lem.\ 7.1]{jantzen96}.
The main distinction is that here the identities are as operators on modules rather than as 
elements of $U\ot U$.

\begin{lem}\label{lem:wef}   Let $\zeta\in Q^+$ and $1\leq i<n$.
Then the following relations hold for operators on tensor products of pairs
of modules in category $\CN$:
\begin{itemize}
\item[(i)] $(\omega_i\ot\omega_i)F_{\zeta}=F_{\zeta}(\omega_i\ot\omega_i)$
      and $(\omega_i'\ot\omega_i')F_{\zeta}=F_{\zeta}(\omega_i'\ot\omega_i')$.
\item[(ii)] $(e_i\ot 1)F_{\zeta} + (\omega_i\ot e_i)F_{\zeta-\alpha_i}=F_{\zeta}
  (e_i\ot 1)+F_{\zeta-\alpha_i}(\omega_i'\ot e_i)$.
\item[(iii)] $(1\ot f_i)F_{\zeta}+(f_i\ot \omega_i')F_{\zeta-\alpha_i}=
     F_{\zeta}(1\ot f_i)+F_{\zeta-\alpha_i}(f_i\ot \omega_i)$.
\end{itemize}
\end{lem}

\begin{proof}
We will verify that (ii) holds.   Identity (iii) is similar to (ii), and (i)
is immediate from (R2) and (R3). In the calculations below, we use Lemma
\ref{lem:xy}(iv)--(vi), (\ref{eqn:lincomb}), and the fact that on a module in category $\CN$  any
element of $U^-_{-\zeta}$ acts as its projection in $\overline{U}^-_{-\zeta}$.
\begin{eqnarray*}
 (e_i\ot 1)F_{\zeta}-F_{\zeta}(e_i\ot 1)&=& \sum_{k=1}^{\overline{d}_{\zeta}}
     (e_iv^{\zeta}_k-v^{\zeta}_ke_i)\ot u^{\zeta}_k\\
   &=& \frac{1}{s-r}\sum_k (p_i'(v^{\zeta}_k)\omega_i'-\omega_ip_i(v^{\zeta}_k))
          \ot u^{\zeta}_k\\
   &=& \frac{1}{s-r}\sum_k\left(\sum_{j=1}^{\overline{d}_{\zeta-\alpha_i}}
        (p_i'(v^{\zeta}_k)\, | \, u^{\zeta-\alpha_i}_j)
      v_j^{\zeta-\alpha_i}\right)\omega_i'\ot u^{\zeta}_k\\
   &&\hspace{.4in}- \ \frac{1}{s-r}\sum_k\omega_i\left(\sum_{j=1}
     ^{\overline{d}_{\zeta-\alpha_i}}(p_i(v^{\zeta}_k)
        \, | \, u_j^{\zeta-\alpha_i})v_j^{\zeta-\alpha_i}\right)\ot u^{\zeta}_k\\
   &=&\sum_k\left(\sum_j(v^{\zeta}_k\, | \, u^{\zeta-\alpha_i}_je_i )v_j^{\zeta-\alpha_i}
     \right)\omega_i'\ot u^{\zeta}_k\\
    &&\hspace{.4in}- \sum_k\omega_i\left(\sum_j(v^{\zeta}_k\, | \, e_iu_j^{\zeta-\alpha_i})
     v_j^{\zeta-\alpha_i}\right)\ot u_k^{\zeta}\\
    &=& \sum_j v_j^{\zeta-\alpha_i}\omega_i'\ot\left(\sum_k(v^{\zeta}_k\, | \,
      u_j^{\zeta-\alpha_i}e_i)u_k^{\zeta}\right)\\
    &&\hspace{.4in} -\sum_j\omega_i v_j^{\zeta-\alpha_i}\ot \left(\sum_k
       (v^{\zeta}_k \, | \, e_iu_j^{\zeta-\alpha_i})u_k^{\zeta}\right)\\
    &=& \sum_{j=1}^{\overline{d}_{\zeta-\alpha_i}}(v_j^{\zeta-\alpha_i}\omega_i'\ot
       u_j^{\zeta-\alpha_i}e_i - \omega_iv^{\zeta-\alpha_i}_j
        \ot e_iu_j^{\zeta-\alpha_i})\\
    &=& F_{\zeta-\alpha_i}(\omega_i'\ot e_i)-(\omega_i\ot e_i)F_{\zeta-\alpha_i}.
\end{eqnarray*}
\end{proof}

Notice that (ii) and (iii) of the above lemma hold even in the cases where
$\overline{U}^+_{\zeta}=0$ but $\overline{U}^+_{\zeta-\alpha_i}\neq 0$,
because in these cases both sides of each equation annihilate modules in category $\CN$.

The following is a modification of \cite[Thm.\ 4.11]{benkart-witherspoon1} or
\cite[Thm.\ 7.3]{jantzen96}. 

\begin{thm}\label{thm:hom} Let $M$ and $M'$ be modules in category $\CN$.
Then the map
$$
  R_{M',M}=F\circ \tilde{f}\circ\tau : M'\ot M\rightarrow M\ot M'
$$
is a homomorphism of $U$-modules.
\end{thm}

\begin{proof}
We must prove that the action of each generator of $U$ commutes with the map $R$.
By Lemma \ref{lem:wef}, $\omega_i$ and $\omega_i'$ commute with $R$.
We will check this for $e_i$ and leave the similar calculation for $f_i$ as an exercise.
We may assume $M=M(\chi)$ and $M'=M'(\psi)$ for some algebra homomorphisms
$\chi,\psi: U^0_{\KK}\rightarrow \KK$. Let $m\in M_{\chi\cdot\hat{\lambda}}$,
$m'\in M'_{\psi\cdot\hat{\mu}}$. By (\ref{eqn:fid}),
\begin{eqnarray*}
  (F\circ\tilde{f}\circ\tau)\Delta(e_i)(m'\ot m)\hspace{-1in}&&\\
     &=& (F\circ \tilde{f})(m\ot e_im'+e_im\ot \omega_im')\\
  &=& f_{\chi,\psi}(\lambda,\mu+\alpha_i) F(m\ot e_im')+f_{\chi,\psi}
     (\lambda+\alpha_i,\mu)F(e_im\ot \omega_im')\\
   &=& f_{\chi,\psi}(\lambda,\mu)\chi(\omega_i')\hat{\lambda}(\omega_i')
           F(1\ot e_i)(m\ot m')\\
  &&\hspace{.5in}+f_{\chi,\psi}(\lambda,\mu)\psi(\omega_i^{-1})
                  \hat{\mu}(\omega_i^{-1})
     F(e_i\ot \omega_i)(m\ot m').
\end{eqnarray*}
Now we may replace $F$ by $\sum_{\zeta\in Q}F_{\zeta}$ or $\sum_{\zeta\in Q}
F_{\zeta-\alpha_i}$, \ $\chi(\omega_i')\hat{\lambda}(\omega_i')m$ by $\omega_i' m$,
and $\psi(\omega_i^{-1})\hat{\mu}(\omega_i^{-1})\omega_im'$ by $m'$.
Thus we obtain the following expression to
which we apply Lemma \ref{lem:wef}(ii):
\begin{eqnarray*}
  f_{\chi,\psi}(\lambda,\mu)\left(\left(\sum_{\zeta}F_{\zeta-\alpha_i}\right)
   (\omega_i'\ot e_i) +\left(\sum_{\zeta}F_{\zeta}\right)
   (e_i\ot 1)\right) (m\ot m')\hspace{-4in} &&\\
&=& f_{\chi,\psi}(\lambda,\mu)\left((e_i\ot 1)\left(\sum_{\zeta}F_{\zeta}\right)
   +(\omega_i\ot e_i)\left(\sum_{\zeta}F_{\zeta-\alpha_i}\right)\right)(m\ot m')\\
  &=& \Delta(e_i) F \tilde{f}\tau(m'\ot m).
\end{eqnarray*}
\end{proof}

We will need the following relation \cite[(5.1)]{benkart-witherspoon1}:
\begin{equation}\label{eqn:yxww}
   (y\omega_{\zeta}' \, | \, x\omega_{\eta})=(y\, | \, x)(\omega_{\zeta}'
    \, | \, \omega_{\eta})
\end{equation}
for all $x\in U^+_{\gamma}$ and $y\in U^-_{-\gamma}$. A similar identity
was proven and used in our Lemma \ref{lem:factor}.

\begin{lem}\label{lem:deltaxy}
If $x\in \overline{U}^+_{\gamma}$ and $y\in \overline{U}^-_{-\gamma}$, then
the following hold as relations of operators on a tensor product of two
modules from category $\CN$:
\begin{itemize}
\item[(i)] $\displaystyle{\Delta(x)=\sum_{0\leq \zeta\leq \gamma} \sum_{i,j}
    (v_i^{\gamma-\zeta}v_j^{\zeta}\, | \, x)u_i^{\gamma-\zeta}\omega_{\zeta}
   \ot u_j^{\zeta}}$,
\item[(ii)] $\displaystyle{\Delta(y)=\sum_{0\leq \zeta\leq\gamma}\sum_{i,j}
    (y\, | \, u_i^{\gamma-\zeta}u_j^{\zeta})v_j^{\zeta}\ot v_i^{\gamma-\zeta}
    \omega_{\zeta}'}$.
\end{itemize}
\end{lem}

\begin{proof}
As in \cite[Lem.\ 5.2]{benkart-witherspoon1}, we may write $\Delta(x)=
\sum_{\zeta,i,j} c_{i,j}^{\zeta} u_i^{\gamma-\zeta}\omega_{\zeta}\ot u_j^{\zeta}$
for some scalars $c^{\zeta}_{i,j}$, this time as an identity of operators on
a tensor product of modules from category $\CN$.
For all $k,l$, and $\nu$, by (\ref{pairprops}) and (\ref{eqn:yxww}), we have
$$
  (v_k^{\gamma-\nu}v_l^{\nu}\, | \, x)=\sum_{\zeta,i,j}c_{i,j}^{\zeta}
   (v_k^{\gamma-\nu} \, | \, u_i^{\gamma-\zeta}\omega_{\zeta})
    (v^{\nu}_l\, | \, u_j^{\zeta})=c^{\nu}_{k,l}.
$$
The proof of (ii) is similar.
\end{proof}

Now let $F^{\op}=\sum_{\gamma\in Q^+}\sum_{i=1}^{\overline{d}_{\gamma}}
u^{\gamma}_i\ot v^{\gamma}_i$, $F^{12}=\sum_{\gamma\in Q^+}\sum_i v_i
^{\gamma}\ot u_i^{\gamma}\ot 1$, similarly $F^{ij}$, and $F^{ij}_f = F^{ij}
\circ \tilde{f_{ij}}$,  where $\tilde{f_{ij}}$ denotes $\tilde f$ applied to the
$i,j$ tensor slots. Recall the notation $R^{ij}$ for $R$ introduced
in Section 2, where the notation $R_{21}$ was reserved for $\tau(R)$.
\medskip

\begin{lem}\label{lem:fop}
The following holds as an identity of operators on a tensor product of three
modules from category $\CN$:
$(\Delta\ot\id)(F^{{\rm{op}}})\circ \tilde{f}_{31}\circ\tilde{f}_{32}
    =F^{31}_f\circ F^{32}_f$.
\end{lem}

\begin{proof} Let $M,M',M''\in \CN$ and assume $M=M(\chi)$, $M'=M'(\psi)$,
and $M''=M''(\phi)$. Let $m\in M_{\chi\cdot \hat{\lambda}}$,
$m'\in M_{\psi\cdot \hat{\mu}}$, $m''\in M''_{\phi\cdot \hat{\nu}}$. The left
side of (i) applied to $m\ot m'\ot m''$ is
\begin{eqnarray*}
  (\Delta\ot\id)(F^{\op})\circ\tilde{f}_{31}\circ\tilde{f}_{32} \hspace{-1in}&&\\
    &=& f_{\phi,\psi}(\nu,\mu)f_{\phi,\chi}(\nu,\lambda)(\Delta\ot\id)
    \left(\sum_{\gamma,k}u_k^{\gamma}\ot v_k^{\gamma}\right)\\
  &=& f_{\phi,\psi}(\nu,\mu)f_{\phi,\chi}(\nu,\lambda)\sum_{\gamma,k}
   \sum_{\zeta,i,j}(v_i^{\gamma-\zeta}v_j^{\zeta}\, | \, u^{\gamma}_k)
     u_i^{\gamma-\zeta}\omega_{\zeta}\ot u_j^{\zeta}\ot v_k^{\gamma}\\
  &=& f_{\phi,\psi}(\nu,\mu)f_{\phi,\chi}(\nu,\lambda)\sum_{\gamma,\zeta,i,j}
    u_i^{\gamma-\zeta}\omega_{\zeta}\ot u_j^{\zeta}\ot\left(\sum_k
    (v_i^{\gamma-\zeta}v_j^{\zeta}\, | \, u_k^{\gamma})v_k^{\gamma}\right)\\
  &=&  f_{\phi,\psi}(\nu,\mu)f_{\phi,\chi}(\nu,\lambda)\sum_{\gamma,\zeta,i,j}
   u_i^{\gamma-\zeta}\omega_{\zeta}\ot u_j^{\zeta}\ot v_i^{\gamma-\zeta}v_j^{\zeta}.
\end{eqnarray*}
On the other hand, by (\ref{eqn:fid}),
\begin{eqnarray*}
  F_f^{31}\circ F_f^{32} (m\ot m'\ot m'')\hspace{-1in}&&\\
   &=& f_{\phi,\psi}(\nu,\mu)\sum_{\eta,\zeta,i,j}f_{\phi,\chi}(\nu-\zeta,\lambda)
    u_i^{\eta}m\ot u^{\zeta}_jm'\ot v_i^{\eta}v_j^{\zeta}m''\\
  &=& f_{\phi,\psi}(\nu,\mu) f_{\phi,\chi}(\nu,\lambda)\sum_{\eta,\zeta,i,j}
   \chi(\omega_{\zeta})\hat{\lambda}(\omega_{\zeta})
    u_i^{\eta}m\ot u_j^{\zeta}m'\ot v_i^{\eta}v_j^{\zeta}m''.
\end{eqnarray*}
Re-summing over $\gamma=\zeta + \eta$, and replacing $\chi(\omega_{\zeta})
\hat{\lambda}(\omega_{\zeta})m$ by $\omega_{\zeta} m$, we
see that (i) holds.
\end{proof}

The next result is the quantum Yang-Baxter equation for $R=F\circ \tilde{f}\circ\tau$.

\begin{thm}\label{thm:QYBE}
Let $M,M',M''\in \CN$. Then $R^{12}\circ R^{23}\circ R^{12}=
R^{23}\circ R^{12}\circ R^{23}$ as maps from $M\ot M'\ot M''$
to $M''\ot M'\ot M$.
\end{thm}

\begin{proof}
If   $\sigma$ is a permutation of $\{1,2,3\}$, let $\tau_{\sigma}$ denote
the corresponding permutation of three tensor factors, that is $\tau_{\sigma}
(m_1\ot m_2\ot m_3)=m_{\sigma^{-1}(1)}\ot m_{\sigma^{-1}(2)}\ot 
m_{\sigma^{-1}(3)}$.  In particular, if $\sigma$ equals the transposition
$(i \, j)$,  we write simply $\tau_{ij}$,  as  $\tau_{ij}$  is just the
same as  $\tau$ applied to tensor slots $i$ and $j$ in that case.  

Note that $\tau_{\sigma} \circ F^{ij}_f=F_f^{\sigma(i)\sigma(j)}\circ
\tau_{\sigma}$ for all $\sigma$, and that $\tilde{f}_{31}\circ
\tilde{f}_{32}\circ F^{12} = F^{12}\circ\tilde{f}_{31}\circ\tilde{f}_{32}$
by two applications of (\ref{eqn:fid}).
We apply these identities, Theorem \ref{thm:hom},  and Lemma \ref{lem:fop}
to obtain the following:
\begin{eqnarray*}
  R^{12}R^{23}R^{12}&=& \tau_{12}\tau_{23}F_f^{31}F_f^{32}R^{12}\\
   &=& \tau_{12}\tau_{23} (\Delta\ot \id)(F^{\op}) \tilde{f}_{31}
   \tilde{f}_{32} F^{12}\tilde{f}_{12}\tau_{12}\\
  &=& \tau_{12}\tau_{23}(\Delta\ot\id)(F^{\op})F^{12}\tilde{f}_{31}
   \tilde{f}_{32} \tilde{f}_{12}\tau_{12}\\
  &=&\tau_{12}\tau_{23}(\Delta\ot\id)(F^{\op}) F^{12}\tilde{f}_{12} \tau_{12}
    \tilde{f}_{32}\tilde{f}_{31}\\
  &=& \tau_{12}\tau_{23} (\Delta\ot\id)(F^{\op}) R^{12} \tilde{f}_{32}
    \tilde{f}_{31}\\
  &=& \tau_{12}\tau_{23} R^{12}(\Delta\ot\id)(F^{\op})\tilde{f}_{32}\tilde{f}_{31}\\
  &=& \tau_{12}\tau_{23} F_f^{12}\tau_{12} F_f^{31}F_f^{32}\\
  &=&F_f^{23} \tau_{12}\tau_{23} \tau_{12} F_f^{31} F_f^{32}\\
  &=& F_f^{23}\tau_{23}\tau_{12}\tau_{23}F_f^{31}F_f^{32}\\
  &=& R^{23}R^{12}R^{23}.
\end{eqnarray*}\end{proof} 
We will need one more lemma in order to obtain the hexagon identities, from
which it will follow that $F$ is a twisting element for $\CN$.

\begin{lem}\label{lem:moreids}
The following are identities of operators on a tensor product of
three modules from $\CN$:
\begin{itemize}
\item[(i)] $(\Delta\ot\id)(F_{\gamma})=\sum_{0\leq\zeta\leq\gamma}
   (F_{\gamma-\zeta})^{23}(F_{\zeta})^{13}(1\ot\omega_{\zeta}'\ot 1)$.
\item[(ii)] $(\id\ot\Delta)(F_{\gamma})=\sum_{0\leq\zeta\leq\gamma}
  (F_{\gamma-\zeta})^{12}(F_{\zeta})^{13} (1\ot\omega_{\zeta}\ot 1)$.
\item[(iii)] $\tilde{f}_{12}\circ (F_{\zeta})^{13} =(F_{\zeta})^{13}
  \circ (1\ot\omega_{\zeta}\ot 1)\circ\tilde{f}_{12}.$
\item[(iv)]  $\tilde{f}_{23}\circ (F_{\zeta})^{13} =(F_{\zeta})^{13}
  \circ (1\ot\omega_{\zeta}'\ot 1)\circ\tilde{f}_{23}.$
\end{itemize}
\end{lem}

\begin{proof}
We will check (ii) and (iv); the proofs of (i) and (iii) are similar.
By Lemma \ref{lem:deltaxy} and (\ref{eqn:lincomb}), considering  operators on
modules we have

\begin{eqnarray*}
  (\id\ot\Delta)(F_{\gamma}) &=& \sum_k\sum_{\zeta,i,j}v_k^{\gamma}
   \ot (v_i^{\gamma-\zeta}v_j^{\zeta}\, | \, u_k^{\gamma}) u_i^{\gamma-\zeta}
      \omega_{\zeta}\ot u^{\zeta}_j\\
   &=& \sum_{\zeta,i,j}\left(\sum_k (v_i^{\gamma-\zeta}v_j^{\zeta}\, | \,
    u^{\gamma}_k)v_k^{\gamma}\right)\ot u_i^{\gamma-\zeta}\omega_{\zeta}
   \ot u_j^{\zeta}\\
   &=& \sum_{\zeta,i,j}v_i^{\gamma-\zeta} v_j^{\zeta}\ot u_i^{\gamma-\zeta}
   \omega_{\zeta}\ot u^{\zeta}_j\\
   &=& \sum_{\zeta} (F_{\gamma -\zeta})^{12}(F_{\zeta})^{13} (1\ot\omega_{\zeta}\ot 1),
\end{eqnarray*}
which proves (ii).

Let $M,M',M''\in \CN$ and $m\in M_{\chi\cdot \hat{\lambda}}$, $m\in M'_{\psi\cdot \hat{\mu}}$,
$m''\in M''_{\phi\cdot \hat{\nu}}$. By (\ref{eqn:fid}),
\begin{eqnarray*}
  \tilde{f}_{23} (F_{\zeta})^{13} (m\ot m'\ot m'') &=& f_{\psi,\phi}(\mu,
    \nu+\zeta)\sum_i v_i^{\zeta}m\ot m'\ot u_i^{\zeta}m''\\
  &=& f_{\psi,\phi}(\mu,\nu)\psi(\omega_{\zeta}')(\omega_{\zeta}'\, | \,
   \omega_{\mu}^{-1})\sum_i v_i^{\zeta}m\ot m'\ot u_i^{\zeta}m''\\
  &=&f_{\psi,\phi}(\mu,\nu) \sum_i v_i^{\zeta}m\ot\omega_{\zeta}'m'\ot u_i^{\zeta}m''\\
  &=& (F_{\zeta})^{13} (1\ot\omega_{\zeta}'\ot 1) \tilde{f}_{23}(m\ot m'\ot m''),
\end{eqnarray*}
which proves (iv).
\end{proof}

Next we will prove the hexagon identities.

\begin{thm}\label{thm:hexagon}
Let $M,M',M''\in \CN$. Then the following are identities of maps from
$M\ot M'\ot M''$ to $M''\ot M\ot M'$ (respectively, $M'\ot M''\ot M$):
\begin{itemize}
\item[(i)] $R^{12}\circ R^{23}=(\id\ot\Delta)(F)\circ \tilde{f}_{12}\circ
  \tilde{f}_{13}\circ \tau_{12}\circ \tau_{23}$.
\item[(ii)] $R^{23}\circ R^{12}=(\Delta\ot \id)(F)\circ \tilde{f}_{23}
  \circ \tilde{f}_{13}\circ\tau_{23}\circ\tau_{12}$.
\end{itemize}
\end{thm}

\begin{proof}
We will prove (ii). The proof of (i) is similar.
Let $m\ot m'\ot m''\in M_{\chi \cdot \hat{\lambda}}\ot M_{\psi \cdot \hat{\mu}}
\ot M_{\phi \cdot \hat{\nu}}$. By Lemma \ref{lem:moreids}, we have
\begin{eqnarray*}
  R^{23}R^{12}(m\ot m'\ot m'') &=& F^{23}\tilde{f}_{23}\tau_{23}F^{12}
  \tilde{f}_{12}\tau_{12}(m\ot m'\ot m'')\\
   &=& F^{23}\tilde{f}_{23}F^{13}\tilde{f}_{13}(m'\ot m''\ot m)\\
  &=& \sum_{\zeta} F^{23}(F_{\zeta})^{13} (1\ot\omega_{\zeta}'\ot 1)
   \tilde{f}_{23}\tilde{f}_{13}(m'\ot m''\ot m)\\
   &=& \sum_{\gamma}\sum_{0\leq\zeta\leq\gamma}(F_{\gamma-\zeta})^{23}
   (F_{\zeta})^{13}(1\ot\omega'_{\zeta}\ot 1)\tilde{f}_{23}\tilde{f}_{13}
   (m'\ot m''\ot m)\\
  &=& \sum_{\gamma}(\Delta\ot \id)(F_{\gamma})\tilde{f}_{23}\tilde{f}_{13}
  (m'\ot m''\ot m)\\
   &=& (\Delta\ot\id)(F)\tilde{f}_{23}\tilde{f}_{13}\tau_{23}\tau_{12}
   (m\ot m'\ot m'').
\end{eqnarray*}
\end{proof}

Finally we show that $F$ is a twisting element for $\CN$.

\begin{thm}\label{deltaf}
Let $M,M',M''$ be modules in category ${\mathcal N}$. Then for
$F$ as defined in (\ref{eqn:F}) we have
$$\bigl[(\Delta\otimes \id)(F)\bigr](F\otimes \id)=
\bigl[(\id \otimes \Delta)(F)\bigr](\id \otimes F)$$
as operators on
$M\ot M'\ot M''$. Thus $F$ is a twisting element for any subcategory
of $\CN$ consisting of $U$-module algebras.
\end{thm}

\begin{proof}
In Theorem \ref{thm:hexagon}, multiply (i) by $R^{12}$ on the right,
multiply (ii) by $R^{23}$ on the right, and apply Theorem \ref{thm:QYBE}
to obtain the identity
$$\bigl[(\Delta\otimes \id)(F)\bigr]\circ\tilde{f}_{23}
\circ\tilde{f}_{13}\circ\tau_{23}\circ\tau_{12}\circ R^{23}=
\bigl[(\id \otimes\Delta)(F)\bigr]\circ \tilde{f}_{12}\circ
\tilde{f}_{13}\circ\tau_{12}\circ
\tau_{23}\circ R^{12}
$$
as functions from $M\otimes M'\otimes M''$ to
$M''\otimes M'\otimes M$, for any $M,M',M''\in \CN$.
It may be checked that
this is equivalent to the identity stated in the theorem,
using (\ref{eqn:fid}). Thus (\ref{twist2}) holds, and (\ref{twist1})
is immediate from the definition of $F$.
\end{proof}

\medskip

We remark that all of the above arguments apply equally well
to $U=U_{r,s}({\mathfrak{sl}}_{n})$,  if  $rs^{-1}$ is
{\em not} a root of unity and  our $U$-module
algebras are (possibly infinite) direct sums of modules arising
from the category $\mathcal{O}$ modules of \cite[Sec.~4]{benkart-witherspoon1}
by shifting weights by some functions $\chi: U^0_{\KK}\rightarrow \KK$.
The function $\tilde{f}$ there merely needs to be replaced by our more
general function $\tilde{f}$ here.
The element $\Theta$  from that paper provides another example of a twisting element
similar to the well-known examples arising from $R$-matrices of
one-parameter quantum groups.

We also note that we may derive Theorem \ref{thm:twist}  as a
special case of Theorem \ref{deltaf}:      When  ${\mathfrak K}=
{\mathbb K}$, the category of modules for $\mathfrak{u}=\mathfrak{u}_{r,s}(\mathfrak{sl}_n)$
is equivalent to a subcategory of $\CN$. The image of
$F$ in $\mathfrak{u}\otimes \mathfrak{u}$ is precisely $R_{f,e}$ by their
definitions.    However, the construction of $R_{f,e}$ in Section
3 shows that it is invertible, while we do not know  if $F$
is invertible in general.

\section{Deformation formulas from infinite quantum groups}

Let $U=U_{r,s}({\mathfrak{sl}}_{n})$, and assume throughout this section that
$r = \theta^y$ and $s = \theta^z$ are both roots of unity with $y,z$ satisfying  
condition (\ref{relprime}).
Let $\KK=\K[[t,t^{-1}]$, Laurent polynomials in $t$
with finitely many negative powers of $t$.
Let $X$ be any finite-dimensional $U$-module such that the extension
$X_{\KK}=X[[t,t^{-1}]$ is a $U_{\KK}$-module in
category $\mathcal N$.
Let $T(X_{\KK})=T_{\KK}(X_{\KK})$ be the tensor algebra of $X_{\KK}$ over $\KK$,
considered as a $U_{\KK}$-module algebra.
As $\mathcal N$ is closed under tensor products and infinite direct sums,
$T(X_{\KK})$ is also in category $\mathcal N$.
Therefore Theorem \ref{deltaf} may be applied to twist the multiplication
of $T(X_{\KK})$.

We wish instead to go further and produce a deformation formula that may be
applied to quotients of $T(X_{\KK})$ to obtain formal deformations.
We define a new action of $U$ on $T(X_{\KK})$ 
denoted $*$, that is $\KK$-linear,
by specifying
\begin{eqnarray*}
  e_i * x &=& (e_i.x)t\\
  f_i * x &=& f_i.x\\
  \omega_i*x &=& (\omega_i.x)t, \qquad \quad \ \omega_i^{-1}*x  =  (\omega_i^{-1}.x)t^{-1} \\
  \omega_i'*x &=& (\omega_i'.x)t, \qquad \quad  (\omega_i')^{-1}*x  =  \big((\omega_i')^{-1}.x\big)t^{-1}
\end{eqnarray*}
for all $x\in X$, and extending to $T(X_{\KK})$
via the $\UK$-module algebra conditions (\ref{modalg1}),
(\ref{modalg2}). (We are just forcing the
generators of $\UK$ to act as scalar multiples of their original actions,
the scalars being powers of $t$.)
The powers of $t$ arising in the new action
$*$ are chosen so that the defining relations of $U$ are
preserved.  
One way to obtain this module structure is to let $U$ act 
on the completion $\K[[t,t^{-1}]\widehat \ot_{\K} X$  of the 
tensor product  in which infinite
sums are allowed.     Here 
$\K[[t,t^{-1}]$ is regarded as a weight module for $U$ with the single weight
$\chi(\omega_i)=\chi(\omega_i')=t$,  so that  $e_i$ and  $f_i$ annihilate
$\K[[t,t^{-1}]$ for each $i$,  and the $U$-action on  
$\K[[t,t^{-1}]\widehat \ot_{\K} X$ is via the coproduct.
Observe that this new $U$-module structure on
$T(X_{\KK})$ yields a module in category
$\mathcal N$, as the only change is in the weights $\chi$.\medbreak

\begin{thm}\label{tv}
The operator $F$ defined in (\ref{eqn:F}) is a universal deformation
formula for any subcategory of $\CN$ (over $\KK=\K[[t,t^{-1}]$)
consisting of $U_{\KK}$-module algebras for which (\ref{eqn:udf}) holds.
\end{thm}

\begin{proof}
This follows immediately from Theorem \ref{deltaf} and the definition
of a universal deformation formula for a category.
\end{proof}

Of course, the $U_{\KK}$-module algebras $T(X_{\KK})$ constructed above
satisfy the hypotheses of the theorem, so by Theorem \ref{udf},
$\mu_{T(X_{\KK})}\circ F$ defines a formal deformation of $T(X)$.

If $X$ is {\em any} finite-dimensional vector space,
the Hochschild cohomology groups ${\rm HH}^i(T(X))$ are $0$
for $i\geq 2$ (e.g.\ see \cite[Prop.\ 9.1.6]{weibel94}).
Thus all formal deformations of $T(X)$ are equivalent to the trivial one, and
so we are interested in applying Theorem \ref{tv} to proper quotients.
For example, we may take the truncated tensor algebra
obtained by letting $p$ be a fixed positive integer and taking the quotient by
the ideal $W_p:=\bigoplus_{k\geq p}X^{\otimes k}$.
The second Hochschild cohomology group of a truncated tensor algebra
is nontrivial. (See for example \cite[Lem.\ 4.1]{cibils}, where we need
to identify our $T(X)/W_p$ with Cibils' $\K \mathcal Q/\mathcal F^p$, $\mathcal Q$ the quiver with
one vertex and $\dim (X)$ loops, and $\mathcal F$ the ideal of the path
algebra $\K \mathcal Q$ generated by the loops.)  
\medbreak

\begin{ex}{\rm  Next we consider some examples that are group crossed products.
Let $V$ be the natural  $n$-dimensional module for $U = U_{r,s}(\mathfrak{sl}_n)$,
defined in Section 3. We will specify an action of $U$ on a
crossed product  $T(V)\# \mathfrak A$,
where $\mathfrak A$ is  the abelian group $(\Z/\ell \Z)^{n-1}$ on generators
$a_1,\ldots,a_{n-1}$, written multiplicatively. Choose arbitrary
$\ell$th roots of unity $\beta_i$ ($1\leq i\leq n-1$), and let an action of
$\mathfrak A$ on $T(V)$ as automorphisms be defined by
$$a_i . v_j = \beta_i v_j$$
for $i=1,\ldots, n-1$ and $j=1,\ldots, n$. The {\em smash product}
$T(V)\# {\mathfrak A}$ is the vector space $T(V)\ot_{\K}
{\mathfrak A}$ with the multiplication
$$ (x\ot a)(y\ot b)=x(a.y)\ot ab,$$
for $x,y \in T(V)$ and $a,b \in \mathfrak A$. Define the following action
of $U$ on the smash product
$\big(T(V)\# \mathfrak A\bigr)[[t,t^{-1}]$:

\begin{eqnarray*}
  e_i.v_j &=& \delta_{i,j-1}v_ia_it,\\
  f_i.v_j&=& \delta_{i,j}v_{i+1},\\
  \omega_i.v_j &=& r^{\delta_{i,j}}s^{\delta_{i,j-1}}
  v_ja_it, \qquad \quad \  \omega_i^{-1}.v_j =  r^{-\delta_{i,j}}s^{-\delta_{i,j-1}}
  v_ja_i^{-1} t^{-1}, \\
  \omega_i'.v_j&=&r^{\delta_{i,j-1}}s^{\delta_{i,j}}
  v_ja_it, \qquad \quad (\omega_i')^{-1}.v_j = r^{-\delta_{i,j-1}}s^{-\delta_{i,j}}
  v_ja_i^{-1} t^{-1},
\end{eqnarray*}
where  $U$ acts trivially on $\mathfrak A$.  
(It is possible to modify this example to involve nontrivial actions
of the $\omega_i$, $\omega_i'$ on $\mathfrak A$.)
It may be checked that $\big(T(V)\# \mathfrak A\bigr)[[t,t^{-1}]$ is a
$U$-module algebra.
Note that $T(V)\# \mathfrak A$ is in category $\mathcal N$, which may be
seen by decomposing $\kc \mathfrak A$ as an algebra  into a direct
sum of copies of $\kc$ (as $\mathfrak A$ is a finite abelian group),
and partitioning the elements of each $V^{\otimes i}\# \mathfrak A$
accordingly,  (or simply by letting $\KK=\kc[[t,t^{-1}]
\mathfrak A$ be the ring of coefficients).
If elements of $\mathfrak A$ are assigned degree 0, it may be verified that $F$
satisfies equation (\ref{eqn:udf}), and so yields a formal
deformation of $T(V)\#{\mathfrak A}$.
This example is similar in some respects to that given in
\cite{caldararu-giaquinto-witherspoon}, although in that
example it was possible to take a further quotient of $T(V)$,
namely a polynomial algebra,  because of the special nature of the
parameters.  

Just as for $T(V)$, formal deformations of $T(V)\# {\mathfrak {A}}$
are necessarily infinitesimally trivial, as the second Hochschild
cohomology group is again trivial.
However, we may again truncate the tensor algebra and consider the
resulting crossed product with $\mathfrak A$ as a $U$-module algebra and
corresponding deformations.} \end{ex}


\begin{thebibliography}{AJS}

\bibitem[AJS]{andersen-jantzen-soergel}H.H.~Andersen, J.C.~Jantzen, and W.~Soergel,
{\it Representations of quantum groups at a $p$th root of unity and of semisimple
groups in characteristic $p$:  Independence of $p$},  Ast\'erisque \textbf{220}
Soci\'et\'e Math. de France  1994.  

\bibitem[BKL]{benkart-kang-lee} G.\ Benkart, S.-J.\ Kang, and K.-H.\ Lee,
{\it On the center of two-parameter quantum groups}, Proc. Roy. Soc. Edinburgh Sect. A,
to appear.

\bibitem[BR]{benkart-roby} G.~Benkart and T. Roby,
{\it Down-up algebras}, J.  
Algebra, {\bf 209} (1998), 305--344; Addendum {\bf 213} (1999), 378.  

\bibitem[BW1] {benkart-witherspoon1}G.\ Benkart and S.\ Witherspoon, {\it Two-parameter
quantum groups and Drinfel'd doubles}, Algebr. Represent.
Theory,  \textbf{7} (2004),   261--286.

\bibitem[BW2]{benkart-witherspoon2} G.\ Benkart and S.\ Witherspoon, {\it Representations
of two-parameter quantum groups and Schur-Weyl duality},
``Hopf Algebras: Proceedings from an International Conference
held at DePaul University'', Bergen, Catoiu, and Chin, eds.,  65--92, Lecture Notes in Pure and Appl. Math. \textbf{237}  Dekker, New York, 2004.

\bibitem[BW3]{benkart-witherspoon3} G.\ Benkart and S.\ Witherspoon, {\it Restricted
two-parameter quantum groups},
{``Representations of Finite Dimensional Algebras
and Related Topics in Lie Theory and Geometry''}, 293--318.
Fields Inst. Commun. \textbf{40},  
Amer. Math. Soc., Providence, RI, 2004.

\bibitem[CGW]{caldararu-giaquinto-witherspoon}A.\ C\u ald\u araru, A.\ Giaquinto,
and S.\ Witherspoon,  {\it Algebraic deformations arising from orbifolds with
discrete torsion},  J.\ Pure Appl.\ Algebra \textbf{187} (2004), 51--70.

\bibitem[CM]{carvalho-musson} P.A.A.B.~Carvalho and I.M.~Musson, {\it Down-up
algebras and their representation theory}, J.\ Algebra \textbf{228} (2000),
286--310. 

\bibitem[C]{cibils} C.\ Cibils, {\it Rigidity of truncated quiver algebras},
Adv.\ Math.\ \textbf{79} (1990), no.\ 1, 18--42.

\bibitem[GZ] {giaquinto-zhang98}A.\ Giaquinto and J.\ Zhang, {\it Bialgebra
actions, twists, and universal deformation formulas,} J.\ Pure Appl.\ Algebra
{\bf 128} (1998), 133--151.

\bibitem[G]{grunspan}C.\ Grunspan, {\it Quantizations of the Witt algebra and
of simple Lie algebras in characteristic $p$},  J.~Algebra \textbf{280} (2004), 145--161.

\bibitem[J]{jantzen96}J.\ C.\ Jantzen, {\it Lectures on Quantum Groups},
vol.\ 6, Graduate Studies in Math., Amer.\ Math.\ Soc., Providence, 1996.

\bibitem[M]{montgomery93}S.\ Montgomery, {\it Hopf Algebras and Their Actions on
Rings}, CBMS Conf.\ Math.\ Publ., vol. \textbf{82}, Amer.\ Math.\ Soc.,
Providence, 1993.

\bibitem[MS]{montgomery-schneider01} S.\ Montgomery and H.J.~Schneider,
{\it Skew derivations of finite-dimensional algebras and actions of the double of the Taft Hopf algebra},   Tsukuba J. Math.  \textbf{25}  (2001),  no. 2, 337--358.


\bibitem[W]{weibel94} C.\ A.\ Weibel, {\it An Introduction to Homological
Algebra}, Cambridge University Press, 1994.

\bibitem[Wi]{witherspoon} S.\ Witherspoon, {\it Skew derivations, Hopf algebras
and deformations of group crossed products}, Commun. Algebra,  to appear,  
({\tt http://www.math.tamu.edu/$\tilde{\hspace{.1cm}}$sjw/pub/skew.pdf}).

\end{thebibliography}
\end{document}